\documentclass{elsart}
\usepackage{amssymb}
\usepackage{ifthen}
\usepackage{hyperref}

\journal{Discrete Mathematics}

\date{}

\newcommand\proofr[1][\hspace{-1em}]{\par{\em Proof\hspace{1em}#1:~}}

\def\proofend{$\blacksquare$\vspace{0.3em}\par}

\makeatletter
\def\@seccntformat#1{\csname the#1\endcsname.\ } 
\makeatother

{\par\addvspace{1mm}{\sc Proof\hspace{1.0ex}{#1}.} }%
{\quad$\blacktriangle$\par\addvspace{1mm}}

\newif\ifNoRemark
\def\addtheorem#1#2#3#4{
\ifthenelse{\equal{#2}{}}{}%
{\ifthenelse{\expandafter\isundefined\csname the#2\endcsname}{\newcounter{#2}}{}}
\newenvironment{#1}[1][\global\NoRemarktrue]
{\par\addvspace{2mm plus 0.5mm minus 0.2mm}\noindent 
\ifthenelse{\equal{#2}{}}{}{\refstepcounter{#2}}%
{\bf #3%
\ifthenelse{\equal{#2}{}}{}{{~\csname the#2\endcsname}}%
 \vphantom{##1}\ifNoRemark.\ \else\ (##1).\fi}\begingroup #4}%
{\endgroup\par\addvspace{1mm plus 0.5mm minus 0.2mm}\global\NoRemarkfalse}
\expandafter\newcommand\csname b#1\endcsname{\begin{#1}}
\expandafter\newcommand\csname e#1\endcsname{\end{#1}}
}

\addtheorem{theorem}{thrm}{Theorem}{\sl}
\addtheorem{lemma}{lmm}{Lemma}{\sl}
\addtheorem{proposition}{prps}{Proposition}{\sl}
\addtheorem{corollary}{crll}{Corollary}{\sl}
\addtheorem{definition}{dfn}{\sl Definition}{}

\providecommand\institute[1]{}

\begin{document}
\begin{frontmatter}
\title{On connection between reducibility of an $n$-ary quasigroup and that of its retracts\thanksref{kropot}}
\author{{Denis~S.~Krotov}},
\ead{krotov@math.nsc.ru}
\author{{Vladimir~N.~Potapov}}
\ead{vpotapov@math.nsc.ru}
\address{Sobolev Institute of Mathematics,\\ prosp. Akademika Koptyuga 4, Novosibirsk, 630090, Russia}
\address{Mechanics and Mathematics Department, Novosibirsk State University,\\ Pirogova 2, Novosibirsk 630090, Russia}

\thanks[kropot]{
This is author's version of the paper in {Discrete Math. 311(1) (2011) 58--66}, DOI:~\href{http://dx.doi.org/10.1016/j.disc.2010.09.023}{10.1016/j.disc.2010.09.023}.
}

\begin{abstract}
An $n$-ary operation $Q:\Sigma^n\to \Sigma$ is called an {$n$-ary
quasigroup} {of order $|\Sigma|$} if in the equation $x_{0}=Q(x_1,
\ldots , x_n)$ knowledge of any $n$ elements of $x_0$, \ldots ,
$x_n$ uniquely specifies the remaining one. An $n$-ary quasigroup $Q$ is (permutably)
reducible if $Q(x_1,\ldots,x_n)=P\left(R(x_{\sigma
(1)},\ldots,x_{\sigma (k)}), \linebreak[1] x_{\sigma
(k+1)},\ldots,x_{\sigma (n)}\right)$ where $P$ and $R$ are
$(n{-}k{+}1)$-ary and $k$-ary quasigroups, $\sigma$ is a permutation,
and $1<k<n$. An $m$-ary quasigroup $R$ is called a retract of $Q$ if
it can be obtained from $Q$ or one of its inverses by fixing $n-m>0$
arguments.

We show that every irreducible $n$-ary quasigroup has an
 irreducible $(n{-}1)$-ary or $(n{-}2)$-ary retract; moreover,
if the order is finite and prime, then it has an
 irreducible $(n{-}1)$-ary retract.
We apply this result to
show that all $n$-ary quasigroups of order $5$ or $7$ whose all
binary retracts are isotopic to $Z_5$ or $Z_7$ are reducible for
$n\geq 4$.
\end{abstract}
\end{frontmatter}

In this paper we prove a reducibility test for $n$-ary quasigroups,
in terms of the reducibility of retracts
(Sections~\ref{s:main}--\ref{s:4-2}), and apply this test to the
classification of the $n$-ary quasigroups of order $5$ or $7$ such
that all their $2$-ary retracts are isotopic to the group $Z_5$ or
$Z_7$ respectively (Section~\ref{s:Z5Z7}).

The test is the following: if all $(n{-}1)$-ary and $(n{-}2)$-ary retracts of an $n$-ary quasigroup $f$ are reducible,
then $f$ is reducible; if all $(n{-}1)$-ary retracts of an $n$-ary quasigroup $f$ of finite prime order are reducible,
then $f$ is reducible. To establish this, we complete the result of the previous paper \cite{Kro:n-3}.
It is worth noting that any test on reducibility makes sense only if there exist irreducible quasigroups.
For different orders and arities, this fact was shown in
\cite{BelSan},\cite{Frenkin},\cite{Borisenko},\cite{Glu76},\cite{Gold75},\cite{Gold76},\cite{AkGo},\cite{KPS:ir}.
The uniqueness of a canonical decomposition of reducible $n$-ary quasigroups
into groups and irreducible quasigroups of smaller arity was proved in \cite{Cher}.

\section{Definitions}\label{s:1}

\begin{definition}
{
 An $n$-ary operation $Q:\Sigma^n\to \Sigma$,
where $\Sigma$ is a nonempty set,
is called
an \emph{$n$-ary qua\-si\-group} or \emph{$n$-qua\-si\-group}
(\emph{of order $|\Sigma|$}) if 
in the equality $x_{0}=Q(x_1, \ldots , x_n)$, knowledge of any $n$ elements
of $x_0$, $x_1$, \ldots , $x_n$ uniquely specifies the remaining one {\rm \cite{Belousov}}.
}
\end{definition}

For the symmetry reasons, we will also indicate an $n$-quasigroup $q$
using the predicate form
$q\langle\cdot\rangle$ defined as
$$
q \langle x_0,x_1,\ldots,x_n\rangle \Leftrightarrow x_{0}=q(x_1,\ldots,x_n).
$$

\begin{definition}
{
 If we assign some fixed values to $l\in \{1,\ldots,n\}$ variables
in the predicate $Q\langle x_0, \ldots , x_n\rangle$
then the $(n-l+1)$-ary predicate obtained corresponds to an $(n-l)$-qua\-si\-group.
Such a qua\-si\-group is called a \emph{retract} or \emph{$(n-l)$-retract} of $Q$.
If $x_0$ is not fixed, the retract is \emph{principal}.
}
\end{definition}

\begin{definition}
{
 By an {\em isotopy} we shall mean
a collection of $n+1$ permutations $\tau_i:\Sigma\to \Sigma$,
$i\in \{0,1,\ldots,n\} $.
$n$-Quasigroups $f$ and $g$ are called {\em isotopic},
if for some isotopy $\bar \tau=(\tau_0,\tau_1,\ldots ,\tau_n)$ we have
$f(x_1,\dots,x_{n}) \equiv \tau_0^{-1} g(\tau_1 x_{1}, \ldots, \tau_{n} x_{n})$, i.\,e.,
$f\langle x_0,x_1,\dots,x_{n}\rangle \equiv g\langle \tau_0 x_{0}, \tau_1 x_{1}, \ldots, \tau_{n} x_{n}\rangle$.
}
\end{definition}

\begin{definition}
{
 An $n$-qua\-si\-group $f$ is termed {\em permutably reducible},
if there exist $m\in\{2,\ldots,n-1\}$,
an $(n-m+1)$-qua\-si\-group $h$,
an $m$-qua\-si\-group  $g$, and a permutation $\sigma: \{1,\ldots,n\} \to \{1,\ldots,n\}$
such that
$$f(x_1,\ldots,x_{n}) \equiv h(g(x_{\sigma(1)},\ldots, x_{\sigma(m)}),
 x_{\sigma(m+1)},\ldots, x_{\sigma(n)}).$$
For short, we will omit the word ``permutably''.
If an $n$-qua\-si\-group is not reducible, then it is {\em irreducible}.
By the definition, all binary
(as well as $1$-ary and $0$-ary) quasigroups are irreducible.
}
\end{definition}

\section{Main results}\label{s:main}
We will prove (Sections~\ref{s:k2}--\ref{s:4-2}) the following two lemmas.
\begin{lemma}\label{l:k2}
Let all the principal $3$- and $4$-retracts of an
$n$-qua\-si\-group $f$ ($n\geq 5$) are reducible. Then
$f$ is reducible.
\end{lemma}

Lemma~\ref{l:k2} was proved in \cite{KroPot:4} for the case of order~4. Indeed,
that proof can be viewed as a simplified version%
\footnote{Even more weak statement whose proof still illustrates 
the proof  of Lemma~\ref{l:k2} is considered in \cite{Kro:swi_sep} for so-called switching separability of graphs,
which is equivalent to the reducibility of the generated quasigroups.}
of the proof of Lemma~\ref{l:k2}.
(In the case of order $4$, all $2$-quasigroups are isotopic to commutative groups;
so, there is no need to work with non-associative $2$-quasigroups and non-commutative groups.
This fact essentially simplifies the proof.)

\begin{lemma}\label{l:n-2}
Let $n\geq 4$.
If an $n$-qua\-si\-group $f$ of finite prime order has an irreducible $(n{-}2)$-retract
and all its $(n{-}1)$-retracts are reducible,
then $f$ is reducible.
\end{lemma}

These lemmas complete the result of \cite{Kro:n-3}, which state that
if the maximum arity $\kappa$ of an irreducible retract
of an $n$-qua\-si\-group belongs to
$\{3,\ldots,n-3\}$,
then the $n$-qua\-si\-group is reducible.
Indeed, Lemma~\ref{l:k2} solves the case $\kappa=2$ and Lemma~\ref{l:n-2}
solves the case $\kappa=n-2$ for prime orders. In summary, we get the following:

\begin{theorem}\label{th}
Let $f$ be an irreducible $n$-qua\-si\-group, $n\geq 4$.
Then $f$ has an irreducible $(n{-}1)$-ary or $(n{-}2)$-ary
retract. Moreover, if the order of $f$ is finite and prime,
then $f$ has an irreducible $(n{-}1)$-retract.
\end{theorem}

As follows from \cite{Kro:n-2}, for every even arity $n\geq 4$ and order $q=4k$
 there exist irreducible $n$-qua\-si\-groups whose entire $(n-1)$-retracts are
reducible.
So, the conditions of Lemma~\ref{l:n-2} and of the last claim of Theorem~\ref{th} cannot be expanded to
all orders.
Nevertheless, the case of odd $n$ and the case of composite order $q \not\equiv 0 \bmod 4$
remain unsolved.

Although Lemma~\ref{l:n-2} does not work for the case of order $4$, this subcase of Theorem~\ref{th}
(proved earlier) was helpful for the characterization of the quasigroups of order $4$ \cite{KroPot:4}.
If it is possible to characterize the class of quasigroups of a larger fixed order or some retract-closed subclass,
the general version of Theorem~\ref{th} can be very helpful, especially for prime orders:
arguing inductively and proving some hypothesis about the structure of irreducible $n$-quasigroups,
we can assume that the $n$-quasigroup has an $(n-1)$-retract that satisfies the hypothesis. The simplest case is when the class consists
of only reducible quasigroups; examples of such classes are shown in
Section~\ref{s:Z5Z7}.

\section{Proof of Lemma~\ref{l:k2}}\label{s:k2}
Let $0$ be some element of $\Sigma$.
Without loss of generality we can consider an $n$-quasigroup $f$ satisfying
$$f(0,\ldots,0,x_i,0,\ldots,0)\equiv x_i$$
for all $i\in\{1,\ldots,n\}$ and $x_i\in\Sigma$ (otherwise, we can apply an appropriate
isotopy).
Assume that all the $3$-ary and $4$-ary retracts of $f$ are reducible.

For binary retracts of $f$, we will use operational symbols like $*$, $\star$, $\diamond$.
For a binary operation $\star$, by $\bar\star$ we denote its conjugate, i.e., $x\star y\equiv y\,{\bar\star}\, x$.
We will write $f|i*(j\star k)$ if $i,j,k\in\{1,\ldots,n\}$
are mutually different and
$$f(0,\ldots,0,x_i,0,\ldots,0,x_j,0,\ldots,0,x_k,0,\ldots,0)\equiv x_i*(x_j\star x_k)$$
(in general, the order of $i,j,k$ can be arbitrary)
where $*$ and $\star$ are binary quasigroups;
if $*$ coincides with $\star$ and is associative (i.\,e., group),
then we will omit the parentheses.
We will write $f|i(jk)$ if $f|i*(j\star k)$ for some $*$ and $\star$.
We will write $f|i[jk]$ if $f|i*(j\star k)$ for some $*$ and $\star$
such that $\star\not\in\{*,\bar *\}$ or $*$ is not associative.
The similar notion will be used with two or four indexes.
Note that
\begin{proposition} \label{p:1234}
For any distinct $i_1, i_2, i_3, i_4 \in \{1,\ldots,n\}$ either
$f|i_a(i_b(i_c i_d))$ or $f|(i_a i_b)(i_c i_d)$ holds, where $\{a,b,c,d\}=\{1,2,3,4\}$.
\end{proposition}

\begin{proposition} \label{p:abc}
{\rm (a)} $f|i(jk)$ if and only if $f|i(kj)$; $f|i[jk]$ if and only if $f|i[kj]$.
\\
{\rm (b)} If $i,j,k\in\{1,\ldots,n\}$
are mutually different then at least one of $f|i(jk)$, $f|j(ik)$, $f|k(ij)$ holds.
\\
{\rm (c)} $f|i[jk]$ if and only if $f|i(jk)$ and neither $f|j(ik)$ nor $f|k(ij)$ holds.
\\
{\rm (d)} $f|i(jk)$ and $f|j(ik)$ mean $f|i*k*j$ for some associative $*$,
which is commutative if and only if $f|k(ij)$.
\end{proposition}

\proofr
(a) is straightforward from the definition. Indeed, $f|i\star(j * k)$ is equivalent to
$f|i\star(k \,{\bar *}\, j)$.

(b) follows from the reducibility of the $3$-retracts of $f$.

(d) Assume $f|i(jk)$ and $f|j(ik)$. Taking into account (a),
we have $f|i*(k \star j)$ and $f|(i \circ k) \diamond j$, i.e.,
\begin{equation}
\label{eq:assoc}
x * (y \star z) \equiv (x \circ y) \diamond z,
\end{equation}
for some $*$, $\star$, $\circ$, and $\diamond$.
(The identity (\ref{eq:assoc}) is known as the general functional equation of associativity and was solved in quasigroups in \cite{Belousov}.)
Substituting consequently $x=0$, $y=0$, and $z=0$,
we find that the operations $*$, $\star$, $\circ$,
and $\diamond$ are the same.
Then, (\ref{eq:assoc}) means the associativity of $*$.
Clearly, the commutativity implies $f|k(ij)$.
Conversely, if $f|k{\scriptstyle\,\triangle\,}(i{\scriptstyle\,\square\,} j)$, then, substituting zero into
$x * y * z \equiv y{\scriptstyle\,\triangle\,}(x {\scriptstyle\,\square\,} z)$, we get $x * y \equiv y{\scriptstyle\,\triangle\,} x$
and $y * z \equiv y{\scriptstyle\,\triangle\,} z$, i.e., $*$ is commutative.

(c) is straightforward from (a), (d), and the definition of the notation $f|i[jk]$.
\proofend

\begin{proposition} \label{p:d}
{\rm (a)} If $f|i[jk]$ and $f|j(kl)$ then $f|i[j(kl)]$ and, in particular,
$i\neq l$, $f|i[jl]$, and $f|i(kl)$.
\\
{\rm (b)} If $f|i*j*l$ and $f|j*k*l$ where $*$ is associative and non-commutative,
then $f|i*j*k*l$.
\end{proposition}
\proofr
(a) Suppose that $f|i[jk]$ and $f|j(kl)$,
i.e., $f|i*(j \star k)$ and $f|j\circ (k \diamond l)$ where
$\star\not\in\{*,\bar *\}$ or $*$ is not associative.
Fixing the $i$th and $l$th variables by zeros,
we see that $\star$ and $\circ$ coincide.
It is easy to check (with the help of Proposition~\ref{p:1234})
that the only variant for the $4$-retract
in the variables $i$, $j$, $k$, $l$ is
$f|i*(j \star (k \diamond l))$. So, we have $f|i*(j \star l)$ and $f|i*(k \diamond l)$,
which prove that $f|i[jl]$ and $f|i(kl)$ (note that we cannot state $f|i[kl]$ because
$*$ can coincide with $\diamond$ or $\bar \diamond$ and be associative).

(b) It is straightforward to check, by Proposition~\ref{p:1234}, that any other variants are inadmissible.
\proofend

\begin{definition}[completely reducible $n$-qua\-si\-groups]
{
An $n$-qua\-si\-group $f$ is called  {\em completely reducible},
if it is reducible and all its principal retracts of arity more than $2$ are reducible.
}
\end{definition}

\begin{lemma}\label{l:phif}
  Assume that all principal $3$- and $4$-retracts of an $n$-qua\-si\-group $f$ are reducible.
Then there exists a completely reducible $n$-qua\-si\-group $\phi_f$ that coincides with $f$
on the $n$-tuples with at most three nonzero elements.
\end{lemma}

\proofr We will argue by induction on $n$. We will show the existence of an `inner' pair,
such that merging this pair enables to reduce $n$.

We say that the pair $\{a,b\}\subset \{1,\ldots,n\}$ is \emph{inner} (\emph{pre-inner}) if
for every $c\in \{1,\ldots,n\}\setminus\{a,b\}$ we have $f|c(ab)$
(respectively, $f|c(ab)$ or $f|a*c*b$ for some associative $*$).

(I). We first prove the existence of a pre-inner pair.

Consider the sequence $b,a_1, \ldots, a_l$ of elements of $\{1,\ldots,n\}$, $1\leq l <n$,
satisfying 
\begin{equation}\label{eq:baa}
f|a_i[a_{i+1} b] \mbox{ for every $i=1,\ldots,l-1$}.
\end{equation}

(*). \emph{We claim that for every $i$ and $j$ such that $1\leq i < j \leq l$
it is true $f|a_i(a_{j} b)$.} Indeed, for $j=i+1$ this is true;
and, by Proposition~\ref{p:d}(a), $f|a_i[a_{i+1} b]$ and $f|a_{i+1}(a_{j} b)$
imply $f|a_i(a_{j} b)$. So, by induction on $j-i$, (*) is proved.

In particular, (*) means that all $a_i$ are pairwise different,
and the length of the sequence is bounded by $n$. Consider a maximal (by length)
sequence ${\mathbf a}:\,b,a_1, \ldots, a_L$ satisfying (\ref{eq:baa}) with $l=L$.

(**). \emph{The pair $b$, $a_L$ is pre-inner.}
Indeed, taking into account Proposition~\ref{p:abc}
and the definition of a pre-inner pair, it is sufficient to avoid
the following two cases:

i) $f|a_L[bc]$ for some $c$; but this contradicts the maximality of $L$,
because $b,a_1, \ldots, a_L, a_{L+1}=c$ also satisfies (\ref{eq:baa}) with $l=L+1$.

ii) $f|b[a_L c]$; in this case we also can show that $L$ is not maximal.
Indeed, consider the sequence ${\mathbf a'}:\,c,a_1, \ldots, a_{L-1}, b, a_L$.
We have \underline{$f|b[a_L c]$}.
If $L=1$, there is nothing to prove.
Otherwise, by Proposition~\ref{p:d}(a), $f|a_{L-1}[a_L b]$ and  $f|b[a_L c]$ yield
\underline{$f|a_{L-1}[b c]$}. 
Inductively, from $f|a_i[a_{i+1} b]$ and $f|a_{i+1}(b c)$ we have
\underline{$f|a_i[a_{i+1} c]$} and $f|a_i(b c)$ (the last is used in
the next induction step) for $i=L-1,L-2,\ldots,2,1$.
So, the underlined expressions confirm that
the sequence ${\mathbf a'}$ satisfies (\ref{eq:baa})
with $l=L+1$, and, thus, $L$ is not maximal.
The contradiction obtained proves (**) and concludes (I).

(II). Now, our goal is to find an inner pair.
We will use the sequence ${\mathbf a}$ considered above.
By the definition, the only possibility for
the pre-inner pair $b$, $a_L$ to be non inner is the existence of  $d_1$ such
that $f|a_L*d_1*b$ where $*$ is the associative non-commutative operation
uniquely specified by $f|a*b$. In that case,

(***) \emph{we state that $f|a_{L-1}[d_1 b]$; i.e., we can replace $a_L$ by $d_1$ in ${\mathbf a}$.}
Indeed, by Proposition~\ref{p:d}(a), $f|a_{L-1}[a_{L} b]$ and $f|a_L(d_1*b)$ imply
$f|a_{L-1}(d_1 b)$; moreover, the operations in
$f|a_{L-1}\star (d_1 * b)$ coincide with that in $f|a_{L-1} \star (a_{L} b)$, and, so,
(***) holds by definition.

So, we get another pre-inner pair $b$, $d_1$. In its turn, it is either inner or there exists $d_2$
such that $f|d_1*d_2*b$. Arguing in such a manner, we construct a sequence
${\mathbf d}:a_L=d_0,d_1,d_2,\ldots$
such that $f|d_{i-1}*d_{i}*b$ and $b$, $d_i$ is a pre-inner pair.
Using Proposition~\ref{p:d}(b), we see by induction that $f|d_j * d_i * b$ for any $j<i$.
This means that all $d_i$ are pairwise different, and the sequence $\mathbf d$ cannot be infinite;
i.e., on some $r$th step we will get an inner pair $b$, $d_r$.

(III) Trivially, the lemma holds for $n = 4$. Assume it holds for $n=k-1$.
Consider the case $n=k$.
As shown above, $f$ has an inner pair of coordinates. Without loss of generality,
it is $\{n-1,n\}$. Let $*$ be the corresponding operation, i.e., $f|(n-1) * n$.
And let the $(n-1)$-qua\-si\-group $q$ be obtained from $f$ by nulling the last argument.
It is straightforward that $\phi_f$ defined as
$\phi_f(x_1,\ldots,x_n)\triangleq\phi_q(x_1,\ldots,x_{n-2},(x_{n-1}*x_n))$
satisfies the condition of the lemma.
\proofend

The last auxiliary statement we need is the following.
\begin{proposition}[{\rm\cite[Theorem~1]{KPS:ir}}] \label{p:24}
 Let $q$ and $g$ be reducible $4$-qua\-si\-groups, and
 let
\begin{eqnarray*}
  q(x,y,z,0)&\equiv& g(x,y,z,0), \\
  q(x,y,0,w)&\equiv& g(x,y,0,w), \\
  q(x,0,z,w)&\equiv& g(x,0,z,w), \ \  \mbox{and}\\
  q(0,y,z,w)&\equiv& g(0,y,z,w).
\end{eqnarray*}
 Then $q\equiv g$.
\end{proposition}

\proofr[of Lemma~\ref{l:k2}]
By Lemma~\ref{l:phif}, there exists a completely reducible $n$-quasigroup $\phi_f$
that coincides with $f$ on the $n$-tuples with at most $3$ non-zero elements.
Using Proposition~\ref{p:24}, it is easy to prove by induction
on the number of non-zero elements in $\bar x
$
that
$\phi_f(\bar x)=f(\bar x)$ for any $\bar x$.

Indeed, let $i,j,k,l$, $1\leq i < j<k<l\leq n$,
be four coordinates with non-zero values in $\bar x$.
Denote by $f_{\bar x; i,j,k,l}$ and $\phi_{\bar x; i,j,k,l}$ the $4$-retracts of $f$ and $\phi_f$,
respectively, with the free coordinates $i,j,k,l$ and the other coordinates fixed by
the values of $\bar x$.
By the induction assumption, $f_{\bar x; i,j,k,l}$ and $\phi_{\bar x; i,j,k,l}$ coincides on the quadruples with at least one zero, i.e.,
satisfy the hypothesis of Proposition~\ref{p:24}. Thus, they are identical;
in particular,
$f_{\bar x; i,j,k,l}(x_i,x_j,x_k,x_l)=\phi_{\bar x; i,j,k,l}(x_i,x_j,x_k,x_l)$,
i.e., $f(\bar x)=\phi_f(\bar x)$.
\proofend

\section{Proof of Lemma~\ref{l:n-2}, case $n>4$}\label{s:n-2}

The line of reasoning of the proof reminds the proof of Lemma~3 in \cite{KroPot:4}.
But there are crucial differences that do not allow to unify these two statements:
in \cite[Lemma~3]{KroPot:4} we deal with $n$-ary quasigroups of order $4$ and prove the
reducibility or semilinearity (some special property of $n$-ary quasigroups of order $4$);
while here we deal with a prime order and prove the reducibility.
To find a special property that would replace the semilinearity in a hypothetical general theorem
is an open problem.

\begin{proposition}[{\rm\cite[Proposition~5.1]{KroPot:4}}] \label{p:canon}
 Assume that a reducible $n$-qua\-si\-group $D$ ($n\geq 3$)
 has an irreducible $(n{-}1)$-retract
 $F\langle x_0,\ldots,x_{n-1}\rangle \equiv  D\langle x_0,\ldots,x_{n-1},0\rangle$.
Then there are $i\in \{0,\ldots ,n\}$ and a $2$-qua\-si\-group $h$
such that $h(x,0)\equiv x$ and
\begin{equation}\label{eq:canon}
D\langle x_0,\ldots,x_{n}\rangle\equiv F\langle x_0,\ldots, x_{i-1},h(x_i,x_{n}), x_{i+1},\ldots, x_{n-1}\rangle.
\end{equation}
\end{proposition}

\begin{proposition} \label{p:aut}
Assume $f$ is an $n$-qua\-si\-group of prime order; $n\geq 3$.
If
$$f\langle x_0,\ldots,x_n \rangle \equiv
f\langle x_0,\ldots,x_{i-1},\mu(x_i),x_{i+1},\ldots,x_{j-1},\nu(x_j),x_{j+1},\ldots,x_n \rangle$$
holds for some different $i,j\in \{0,\ldots,n\}$ and some pair of non-identity permutations $(\mu,\nu)$,
then $f$ is reducible.
\end{proposition}
\proofr
 Without loss of generality assume that $i=1$, $j=2$.
Put
\begin{eqnarray*}
\alpha(x,y) &\triangleq& f(x,y,\bar 0), 
\\
\beta(x,\bar z) &\triangleq& f(x,0,\bar z),\qquad \bar z\triangleq(z_1,\ldots ,z_{n-2}),\\
\gamma(x) &\triangleq& f(x,0,\bar 0). 
\end{eqnarray*}
Note that
$$\gamma^{-1}(\alpha(x,0))\equiv x.$$
Assume that the pair $(\mu, \nu)$ satisfies the hypothesis of Proposition.
Then $\alpha(x,y) \equiv \alpha(\mu x,\nu y)$, and $\alpha(x,y) \equiv \alpha(\mu^s x,\nu^s y)$
for all natural $s$.
By the definition of quasigroup, $x=\mu^s x$ if and only if $y=\nu^s y$.
So, it is easy to see
that the permutations $\mu$ and $\nu$ consist of cycles of the same length;
since the order is prime, $\mu$ and $\nu$ are cyclic permutations.
As corollary, we have the following:

(*) For each
$v\in\Sigma$ there exist permutations $\rho_v,\tau_v:\Sigma\to \Sigma$
such that $f(x,y,\bar z)\equiv f(\rho_v x,\tau_v y,\bar z)$ and $\tau_v v=0$
(in other words, the group of permutations $\tau$ admitting
$f(x,y,\bar z)\equiv f(\rho x,\tau y,\bar z)$
for some $\rho$ acts transitively on $\Sigma$, i.e., has only one orbit).

 Then,
\begin{eqnarray*}
f(x,y,\bar z)
&\equiv&
f(\rho_y x,\tau_y y,\bar z)
\equiv
f(\rho_y x,0,\bar z)
\equiv
\beta(\rho_y x,\bar z)
\\&\equiv&
\beta(\gamma^{-1}\alpha(\rho_y x,0),\bar z)
\equiv
\beta(\gamma^{-1}\alpha(\rho_y x,\tau_y y),\bar z)
\equiv
\beta(\gamma^{-1}\alpha(x,y),\bar z)
\end{eqnarray*}
and thus $f$ is reducible provided $n\geq 3$.
\proofend

{\it Proof of Lemma~\ref{l:n-2}, case $n>4$: }
Assume $C$ is an $n$-qua\-si\-group.
Assume all the $(n-1)$-retracts of $C$ are reducible
and $C$ has an irreducible $(n-2)$-retract $E$.
Without loss of generality assume that
$$ E\langle x_0,\ldots,x_{n-2} \rangle \equiv C \langle x_0,\ldots,x_{n-2},0,0 \rangle. $$
We will use the following notation for retracts of $C$:

\begin{eqnarray}
  E_{a,b}\langle x_0,\ldots,x_{n-2} \rangle & \triangleq & C \langle x_0,\ldots,x_{n-2},a,b \rangle, \label{eq:E}\\
  A_b\langle x_0,\ldots,x_{n-2},y \rangle & \triangleq & C\langle x_0,\ldots,x_{n-2},y,b \rangle, \label{eq:A}\\
  B_a\langle x_0,\ldots,x_{n-2},z \rangle & \triangleq & C\langle x_0,\ldots,x_{n-2},a,z \rangle. \label{eq:B}
\end{eqnarray}

(*) Since $A_0$ is reducible and fixing $y:=0$ leads to the irreducible $E$,
by Proposition~\ref{p:canon} we have
\begin{equation}\label{eq:A0}
A_0\langle x_0,\ldots,x_{n-2},y \rangle
\equiv
E\langle x_0,\ldots,x_{i-1},h(x_i,y),x_{i+1},\ldots,x_{n-2} \rangle
\end{equation}
for some $i\in \{0,\ldots n-2\} $ and $2$-qua\-si\-group $h$ such that $h(x_i,0)\equiv x_i$.

From (\ref{eq:A0}), we see that all the retracts
$E_{a,0}$, $a\in\Sigma$, are isotopic to $E$. Similarly, we can get the following:

(**) All the retracts $E_{a,b}$, $a,b\in\Sigma$ are isotopic to $E$.

Then, similarly to (*), we have that for each $b\in\Sigma$
$$A_b \langle x_0,\ldots,x_{n-2},y \rangle
\equiv
E_{0,b} \langle x_0,\ldots,x_{i_b-1},h_b(x_{i_b},y),x_{i_b+1},\ldots,x_{n-2} \rangle
$$
for some $i\in \{0,\ldots n-2\} $ and $2$-qua\-si\-group $h_b$ such that $h_b(x,0)\equiv x$.

(***) We claim that $i_b$ does not depend on $b$.
Indeed, assume, for example, that $i_1=0$ and $i_2=1$, i.e.,
\begin{eqnarray*}
A_1 \langle x_0,\ldots,x_{n-2},y \rangle
&\equiv&
E_{0,1} \langle h_1(x_0,y),x_1,x_2,\ldots,x_{n-2} \rangle,\\
A_2 \langle x_0,\ldots,x_{n-2},y \rangle
&\equiv&
E_{0,2} \langle x_0,h_2(x_1,y),x_2,\ldots,x_{n-2} \rangle.
\end{eqnarray*}
Then, fixing $x_0$ in the first case leads to a retract isotopic to $E$;
fixing $x_0$ in the second case leads to a reducible retract (recall that $n\geq 5$).
But, analogously to (**), these two retracts are isotopic;
this contradicts the irreducibility of $E$ and proves (***).

Without loss of generality we can assume that $i_b=0$, i.\,e.,
\begin{equation}\label{eq:A0E}
A_b \langle x_0,x_1,\tilde x_2,y \rangle \equiv E_{0,b} \langle
h_b(x_0,y),x_1,\tilde x_2 \rangle;
\qquad\qquad\mbox{here and later $\tilde x_2\triangleq (x_2,\ldots,x_{n-2})$}.
\end{equation}
Similarly, we can assume without loss of generality that either
\begin{equation}\label{eq:B0E}
B_a \langle x_0,x_1,\tilde x_2,z \rangle \equiv E_{a,0} \langle
g_a(x_0,z),x_1,\tilde x_2 \rangle
\end{equation}
or
\begin{equation}\label{eq:B1E}
B_a \langle x_0,x_1,\tilde x_2,z \rangle \equiv E_{a,0} \langle
x_0,g_a(x_1,z),\tilde x_2 \rangle
\end{equation}
where $2$-qua\-si\-groups $g_a$ satisfy $g_a(x,0)\equiv x$.

Assuming that (\ref{eq:B0E}) holds, we derive
\begin{eqnarray}
C \langle x_0,x_1,\tilde x_2,y,z\rangle
&\stackrel{(\ref{eq:A})}{\equiv}& \nonumber
A_z\langle x_0,x_1,\tilde x_2,y \rangle
\\
  &\stackrel{(\ref{eq:A0E})}{\equiv}& \nonumber
  E_{0,z}\langle h_z(x_0,y),x_1,\tilde x_2\rangle
\\
  &\stackrel{(\ref{eq:E}),(\ref{eq:B})}{\equiv}& \label{eq:4-5}
  B_{0}\langle h_z(x_0,y),x_1,\tilde x_2, z\rangle
\\
  &\stackrel{(\ref{eq:B0E})}{\equiv}& \nonumber
  E_{0,0} \langle g_0(h_z(x_0,y),z),x_1,\tilde x_2 \rangle,
\end{eqnarray}
which means that $C$ is reducible, because
$f(x,y,z) \triangleq g_0(h_z(x,y),z)$ must be a $3$-qua\-si\-group
(fixing $x_1$ and $\tilde x_2$, we see that it is a retract of $C$).
So, it remains to consider the case (\ref{eq:B1E}).
Consider two subcases:

Case 1: $g_a$ does not depend on $a$; denote $g \triangleq g_0=g_a$ for all $a\in\Sigma$.
Then, repeating the three steps of (\ref{eq:4-5}),
we derive
$$
C\langle x_0,x_1,\tilde x_2,y,z \rangle
  \stackrel{(\ref{eq:B})}{\equiv}
B_y\langle x_0,x_1,\tilde x_2,z \rangle
  \stackrel{(\ref{eq:B1E})}{\equiv}
E_{y,0} \langle x_0,g(x_1,z),\tilde x_2,z \rangle
  \stackrel{(\ref{eq:E}),(\ref{eq:A})}{\equiv}
A_0 \langle x_0,g(x_1,z),\tilde x_2,y \rangle
$$
and see that $C$ is reducible.

Case 2: for some fixed $a$ we have $g_0\neq g_a$;
denote $s_i(x)\triangleq g_0(x,i)$, $t_i(x)\triangleq g_a(x,i)$,
and $r_i(x) \triangleq h_i(x,a)$.
From (\ref{eq:A0E}), we see that
\begin{eqnarray}
E_{a,0}\langle x_0,x_1,\tilde x_2 \rangle
&\equiv & \label{eq:E-1}
E_{0,0}\langle r_0(x_0),x_1,\tilde x_2 \rangle,
\\
E_{a,b}\langle x_0,x_1,\tilde x_2 \rangle
&\equiv & \label{eq:E-2}
E_{0,b}\langle r_b(x_0),x_1,\tilde x_2 \rangle.
\end{eqnarray}
From (\ref{eq:B1E}), we see that
\begin{eqnarray}
E_{0,b}\langle x_0,x_1,\tilde x_2 \rangle
&\equiv & \label{eq:E-3}
E_{0,0}\langle x_0,s_b(x_1),\tilde x_2 \rangle,
\\
E_{a,b}\langle x_0,x_1,\tilde x_2 \rangle
&\equiv & \label{eq:E-4}
E_{a,0}\langle x_0,t_b(x_1),\tilde x_2 \rangle.
\end{eqnarray}
Applying consequently (\ref{eq:E-3}), (\ref{eq:E-2}), (\ref{eq:E-4}), and (\ref{eq:E-1}),
we find that for each $b$ the retract  $E=E_{0,0}$ satisfies
$$
E\langle x_0,x_1,\tilde x_2 \rangle
\equiv E_{0,b}(\ldots)
\equiv E_{a,b}(\ldots)
\equiv E_{a,0}(\ldots)
\equiv E\langle r_b^{-1} r_0 x_0,s_b^{-1} t_b x_1,\tilde x_2 \rangle.
$$
By Proposition~\ref{p:aut}, the irreducibility of $E$ means that
$s_b^{-1} t_b=Id$
for every $b$. But this contradicts the assumption $g_0\neq g_a$ and, so, proves
that Case 2 is not admissible.
\proofend

\section{Proof of Lemma~\ref{l:n-2}, case $n=4$}\label{s:4-2}

Let $\Sigma=\{0,1,\ldots,k-1\}$. In this section, an important role is played by the $2$-qua\-si\-group $x+y\bmod k$,
which will be denoted by $Z_k$.
All arithmetic operations (addition, multiplication)
will be performed modulo $k$ where $k$ is the order of considered quasigroups.

\begin{proposition}\label{prozk1}
   Assume $f$ is a $2$-quasigroup of a prime order $k$.
   If there are non-identity permutations
   $\nu, \mu:\Sigma\to\Sigma$ such that either $f(\mu (x),\nu(y))\equiv f(x,y)$, or
$f(\mu (x), y)\equiv \nu (f(x,y))$, or
$f( x,\nu (y))\equiv \mu (f(x,y))$,
then $f$ if isotopic to $Z_k$.
\end{proposition}
\proofr
Consider only the case $f(\mu (x),\nu (y))\equiv f(x,y)$, because the other cases can be reduced to it by considering $f^{(2)}$ or $f^{(1)}$ respectively, where $f^{(i)}$ is the inversion of $f$ in the $i$th argument.

We first note that the permutations $\mu$ and $\nu$ must be cycles from $k$
elements (a similar statement can be found in the proof of Proposition~\ref{p:aut}).

Now, let us prove that $f$ is isotopic to $Z_k$.
There are permutations $\alpha$ and $\beta$ such that
$\alpha\mu\alpha^{-1} = \epsilon$ and $\beta\nu\beta^{-1} = \epsilon^{-1}$,
where $\epsilon$ corresponds to the addition of $1$ $\bmod\ k$
(explicitly, $\alpha$ can be recursively defined by
$\alpha(0)\triangleq 0$, $\alpha(\mu (x))\triangleq \alpha (x) +1$; similarly for $\beta$).
Then the quasigroup $g(x,y) \equiv f(\alpha^{-1} (x), \beta^{-1} (y))$ meets
$g(\epsilon (x),\epsilon^{-1}(y)) \equiv g(x,y)$. Next, define
$h(x,y) \equiv \gamma^{-1}(g( x,  y))$
where $\gamma(x) \triangleq g(x,0)$. Let us check now that we have got $Z_k$:
$
h(x,y)\equiv h(\epsilon^y (x), \epsilon^{-y} (y))\equiv h(x+y,0)\equiv x+y.
$
 \proofend

\begin{proposition}\label{prozk21}
   Let $|\Sigma|$ be prime.
   Assume $\tau (x + a)\equiv \tau (x) +b$ for some permutation $\tau:\Sigma\to\Sigma$
   and $a,b\in \Sigma$. Then either $\tau (x) \equiv c \cdot x +d$
   for some  $c,d\in \Sigma$ or $a=b=0$.  \end{proposition}
\proofr If $a=0$ then obviously $b=0$. Let $a\neq 0$.
From $\tau (x + a)\equiv \tau (x) +b$ we get, by induction, that
$\tau (y +
i\cdot a)\equiv \tau(y)+i\cdot b$ for every $i\in\Sigma$.
Further,
$$\tau (x)\equiv \tau(0+(x/a)\cdot a) \equiv  \tau (0) + (x/a)\cdot b \equiv (b/a)\cdot x + \tau (0).$$
 \proofend

\begin{proposition}\label{prozk3} Let a quasigroup $f$ be isotopic to $Z_k$.
If $f(0,x)\equiv f(x,0)\equiv x$ and $f(x,1)\equiv x+1$, then $f(x,y)\equiv x+y$.
\end{proposition}

\proofr We have
$\alpha (f(x,y))\equiv
\beta (x) +\gamma (y)$
for permutations $\alpha,\beta,\gamma:\Sigma\to\Sigma$.
From $f(0,y)\equiv y$ we get
\begin{equation} \label{eq:g=a-b}
\gamma(y) \equiv \alpha(y) - \beta( 0).
\end{equation}
From $f(x,0)\equiv x$ we get
\begin{equation} \label{eq:b=a-g}
\beta (x) \equiv \alpha (x) - \gamma (0).
\end{equation}
From $f(x,1)\equiv x+1$ we get
$\alpha(x+1) \equiv \beta (x)+\gamma (1) \equiv
\alpha (x) - \gamma (0)+\gamma (1)$.
Hence, by induction,
$\alpha (z)\equiv z\cdot (\gamma(1)-\gamma(0))+\alpha(0)$.
Consequently,
\begin{equation} \label{eq:a=a+a-a}
\alpha (x+y)\equiv\alpha (x)+\alpha (y)-\alpha(0).
\end{equation}
Using (\ref{eq:g=a-b})--(\ref{eq:a=a+a-a}), we derive
$$\alpha (f(x,y)) \equiv
\beta(x)+\gamma(y) \equiv
\alpha (x) - \gamma (0)+ \alpha(y) - \beta( 0) \equiv
\alpha (x) +\alpha(y) -\alpha( 0) \equiv \alpha(x+y),
$$
which proves that $f(x,y)\equiv x+y$.
 \proofend

\begin{corollary} \label{cor:z5z5} Let $k$ be a prime integer;
let $2$-quasigroups $f$ and $g$ be isotopic to $Z_k$;
and let for some fixed $a$, $a'$, $b$, $b'$, $a\neq b$, the identities
$f(x,a)\equiv g(x,a')$ and $f(x,b)\equiv g(x,b')$ hold.
Then $f(x,y)\equiv g(x,\pi(y))$ for some permutation $\pi$.
\end{corollary}
\proofr
Proposition~\ref{prozk3} is a partial case of the considered statement;
we will reduce the general case to this partial case.
Note that the statement is equivalent to the similar statement
for quasigroups the $f'$ and $g'$ obtained from $f$ and $g$
by applying some (common) isotopy. Using this observation,
we may

1) assume without loss of generality that $g(x,y)\equiv x+y$;

2) using an isotopy $(\mathrm{Id},\alpha,\alpha^{-1})$ where
$\alpha$ is the addition of a constant, assume that $a'=0$;

3) using an isotopy $(\beta,\beta,\beta)$
where $\beta$ is the multiplication by a constant, assume that  $b'=1$
(here we use that $k$ is prime).

Furthermore, we may replace $f$ by $f''$ where $f''(x,y) \equiv f(x,\gamma(y))$ for some permutation $\gamma$.
This permutation can be chosen in such a way that $f''(0,y) \equiv y$
(since $a'=0$ and $b'=1$, we will also get $a=0$, $b=1$).
Then the statement follows from Proposition~\ref{prozk3}.
\proofend

We leave the proof of the following simple proposition; a similar statement can be found in \cite[Proposition~2]{Kro:n-3}.

\begin{proposition}\label{proraz4} Let $f$ be a $4$-quasigroup of finite order,
and let all its $3$-retracts are reducible. Then all  $2$-retracts of $f$
that obtained by fixing the same arguments with different values are mutually isotopic.\end{proposition}

\begin{proposition}\label{proraz6} Let $f$ be a $4$-quasigroup of finite order $k$
and let there be $2$-quasigroups $q$, $r$, $g_z$, $h_z$, $z\in\Sigma$,
such that
$$ f(x_1,x_2,x_3,z)\equiv g_z(x_3,h_z(x_1,x_2))\qquad \mbox{and}\qquad f(0,x_2,x_3,z)\equiv q(x_3,r(z,x_2)).$$
Then $f$ is reducible. \end{proposition}

\proofr
Setting $x_1=0$, we get
$g_z(x_3,h_z(0,x_2))\equiv q(x_3,r(z,x_2))$.
Hence for every $z$ we have $g_z(x_3,x_2)\equiv q(x_3,\pi_z(x_2))$ for some permutation $\pi_z$.
Then $f(x_1,x_2,x_3,z)\equiv q(x_3,\pi_z^{-1}(h_z(x_1,x_2)))$, which implies that
$f$ is reducible.
\proofend

\proofr[of Lemma~\ref{l:n-2}, case $n=4$]
Note that  any $4$-quasigroup of order $2$ or $3$ is isotopic to the iterated group $Z_k$ and, so, reducible.
Suppose $k>3$.
Let $f:\Sigma^4\to\Sigma$ be a $4$-quasigroup and let all its
$3$-retracts be reducible.
For every $i=1,2,3,4$ and $a\in \Sigma$ define
\def\uot{\mathrm{out}}
$\uot(i,b)$ as the order number of the argument of $f$ that will be \emph{outer} in the decomposition of $f(x_1,x_2,x_3,x_4)$ with $x_i$ fixed by $b$. E.g., $\uot(1,2)=3$ means that
$f(2,x_2,x_3,x_4)\equiv g(x_3,h(x_2,x_4))$ for some $g$ and $h$.
If the corresponding retract admits two or three decompositions with different outer arguments, we,
for definiteness, require $\uot(i,b)$
to be  minimum  among the admissible values.

By simple counting arguments, it is easy to find that one
of the two following statements is true:

(1) There are two different positions $i,j\in\{1,2,3,4\}$ and values $a,c,d\in\Sigma$, $c\neq d$, such that
$\uot(i,a)=\uot(j,c)=\uot(j,d)$.

(2) The function $\uot(i) \triangleq \uot(i,b)$ does not depend on $b$ and bijective with respect to $i$.
Moreover, we can assume that $\uot(i)$ is a cyclic permutation,
 because otherwise the inversion of $f$ in any position (say, the first) satisfies (1).

Suppose (1) holds.
Consider, without loss of generality, the following decompositions ($j=3$, $i=4$, $a=c=0$, $d=1$):
\begin{eqnarray}
f(x,y,0,w)&\equiv& g_0(y,h_0(x,w)),\label{eq:f0}\\
f(x,y,1,w)&\equiv& g_1(y,h_1(x,w)),\label{eq:f1}\\
f(x,y,z,0)&\equiv& g(y,h(x,z))\label{eq:fx}.
\end{eqnarray}
Consider the retracts $\varphi_{z,w}(x,y)\triangleq
f(x,y,z,w)$. The further arguments will be divided into two subcases.

(1a) \emph{The quasigroup $\varphi_{0,0}$ is isotopic to $Z_k$}.

Setting $z=0$ and $w=0$ in (\ref{eq:f0}) and (\ref{eq:fx}),
we see that $g_0$ and $g$ are isotopic to $Z_k$;
moreover, $g_0(y,t)\equiv g(y,\tau(t))$ for some permutation $\tau$.
The similar is true for $g_1$. So, we can assume without loss of generality that
 (\ref{eq:f0}) and (\ref{eq:f1}) hold with
\begin{equation}\label{eq:g0g1g}
g_0=g_1=g.
\end{equation}

Now consider the $2$-retracts $q_{x,w}(y,z)\equiv f(x,y,z,w)$.
From (\ref{eq:fx}) we see that if $w=0$ then $q_{x,w}$ is isotopic to $g$;
hence, it is isotopic to $Z_k$. By Proposition~\ref{proraz4}, it is true for every
$x$ and $w$.
Next, from (\ref{eq:f0}) and (\ref{eq:g0g1g}) we have that for fixed $x$ and $w$
$q_{x,w}(y,0)\equiv g(y,H_{0,x,w})$ for some $H_{0,x,w}$ (explicitly, $H_{0,x,w}\triangleq h_0(x,w)$).
Similarly,
$q_{x,w}(y,1)\equiv g(y,H_{1,x,w})$ (with $H_{1,x,w}\triangleq h_1(x,w)$).
By Corollary~\ref{cor:z5z5} we have $q_{x,w}(y,z)\equiv g(y,\pi_{x,w}(z))$
for some permutation $\pi_{x,w}$. But this means  $f(x,y,z,w)\equiv g(y,\pi(x,z,w))$, where
$\pi(x,z,w) \triangleq \pi_{x,w}(x)$, which proves the reducibility of $f$.

(1b) \emph{The quasigroup $\varphi_{0,0}$ is not isotopic to $Z_k$}.
Let us show that

(*) \emph{all the retracts of $f$ obtained by fixing the fourth argument have the same (second) outer argument}.
Seeking a contradiction, consider, as an example, the alternative decomposition
\begin{equation}
f(x,y,z,2)\equiv g'(x,h'(y,z))\label{eq:fx'}.
\end{equation}

From (\ref{eq:f0}) we get
$\varphi_{0,2}(x,y) \equiv \varphi_{0,0}(\pi_0^{-1}(\pi_2(x)),y)$,
where $\pi_i(x) \equiv h_0(x,i)$ (we do not need this explicit form below).
Similarly, from (\ref{eq:fx'}),
$\varphi_{1,2}(x,y) \equiv \varphi_{0,2}(x,\nu_0^{-1}(\nu_1(y)))$,
where $\nu_i(y) \triangleq h'(y,i)$ (observe that $\nu_0\neq\nu_1$).
Similarly,
we express $\varphi_{1,0}$ in terms of $\varphi_{1,2}$, using (\ref{eq:f1});
and, finally, $\varphi_{0,0}$ in terms of $\varphi_{1,0}$ (\ref{eq:fx}).
In summary, we obtain
$\varphi_{0,0}(x,y)\equiv\varphi_{0,0}(\mu (x),\nu (y))$ for some permutation $\mu$ and
 $\nu = \nu_0^{-1} \nu_1 \neq\mathrm{Id}$.
Using Proposition~\ref{prozk1}, we conclude that $\varphi_{0,0}$ is isotopic to $Z_k$,
which contradicts the hypothesis of the considered subcase.
Similarly, the other way to decompose $f(x,y,z,2)$ (with outer $x$) also
leads to a contradiction; moreover, $2$ can be replaced by any nonzero value.
Claim (*) is proved.

Now, the reducibility of $f$ is straightforward from
Proposition~\ref{proraz6}.

Let (2) hold. Without loss of generality let $\uot(3)=2$,
$\uot(2)=4$, $\uot(4)=1$, $\uot(1)=3$. I.e.,
\begin{eqnarray}
f(x,y,z,w)
 &\equiv& g_z(h_z(x,w),y)  \label{eq:Z} \\
 &\equiv& q_w(x,r_w(z,y))  \label{eq:W} \\
 &\equiv& \alpha_x(z,\beta_x(w,y))  \label{eq:X} \\
 &\equiv& \gamma_y(\delta_y(x,z),w). \label{eq:Y}
\end{eqnarray}
for some quasigroups $g_t$, $h_t$, $q_t$, $r_t$, $\alpha_t$,
$\beta_t$, $\gamma_t$, $\delta_t$, $t\in \Sigma$.
Without loss of generality we may assume that
$f(0,0,0,v) \equiv f(0,0,v,0) \equiv  f(0,v,0,0) \equiv f(v,0,0,0) \equiv v$.
Moreover, the quasigroups
$g_0$, $h_0$, $q_0$, $r_0$, $\alpha_0$, $\beta_0$,
$\gamma_0$, $\delta_0$
can be chosen in such a way that every $p$ from them satisfies
$p(0,v)\equiv p(v,0) \equiv v$ (see, e.g., \cite[Proposition~12]{PotKro:asymp}).
Fixing two variables by zeros in
(\ref{eq:Z})--(\ref{eq:Y}), we derive that
$g_0 = q_0 = \beta_0=\delta_0$ and
$h_0 = r_0 = \alpha_0=\gamma_0$. Further arguments will differ
depending on whether all these quasigroups are isotopic to $Z_k$ or not.

Again, we consider the retracts
$\varphi_{z,w}(x,y)\triangleq f(x,y,z,w)$ and divide the arguments
into two subcases.

(2a)
{\em The quasigroups $\varphi_{z,w}$ (as well as $g_0$, $q_0$, $\beta_0$, $\delta_0$) are not isotopic to $Z_k$.}
Consider the retracts $\varphi_{0,0}$, $\varphi_{0,b}$,
$\varphi_{a,0}$, $\varphi_{a,b}$ for some  $a$ and $b$ from $\Sigma\setminus\{0\}$.
From (\ref{eq:Z})  we have
\begin{eqnarray*}
\varphi_{0,0}(x,y)&\equiv& f(x,y,0,0) \equiv
g_0(h_0(x,0),y),
\\
\varphi_{0,b}(x,y)&\equiv& f(x,y,0,b) \equiv
g_0(h_0(x,b),y).
\end{eqnarray*}
So,
\begin{equation}\label{eq:000b}
\varphi_{0,b}(x,y) \equiv \varphi_{0,0}(\sigma_{0,b}(x),y)
\end{equation}
for some permutation $\sigma_{0,b}$ (explicitly,
$\sigma_{0,b}(x) \triangleq h_0(h_0^{-1}(x,0),b)$\,).
Similarly, from (\ref{eq:Z}) and (\ref{eq:W}) we derive
\begin{eqnarray}
\label{eq:a0ab}
\varphi_{a,b}(x,y) \equiv \varphi_{a,0}(\sigma_{a,b}(x),y) \\
\label{eq:00a0}
\varphi_{a,0}(x,y) \equiv \varphi_{0,0}(x,\tau_{a,0}(y)) \\
\label{eq:0bab}
\varphi_{a,b}(x,y) \equiv \varphi_{0,b}(x,\tau_{a,b}(y))
\end{eqnarray}
for some $\sigma_{a,b}$, $\tau_{a,0}$, $\tau_{a,b}$.
We will show that $\sigma_{0,b}=\sigma_{a,b}$
and $\tau_{a,0}=\tau_{a,b}$. Indeed, from (\ref{eq:000b})--(\ref{eq:0bab}) we have
$$
\varphi_{0,0}(x,y)\equiv
\varphi_{0,0}(\sigma_{0,b}^{-1}(\sigma_{a,b}(x)), \tau_{a,b}(\tau^{-1}_{a,0}(y))),
$$
and, by Proposition~\ref{prozk1},
taking into account the condition of subcase (2a),
we conclude that
$\sigma_{0,b}^{-1}\sigma_{a,b}
=\tau_{a,b}\tau^{-1}_{a,0}=\mathrm{Id}$.

Now, we have
$$
f(x,y,z,w) \equiv \varphi_{z,w}(x,y)
\equiv \varphi_{0,0}(\sigma_{0,w} (x),\tau_{z,0} (y)),
$$
which implies that $f$ is reducible.

(2b) {\em The quasigroups $\varphi_{z,w}$ are isotopic to $Z_k$}.
Assume without loss of generality that $g_0$, $q_0$, $\beta_0$, $\delta_0$ coincide with $Z_k$.
We may also assume that the quasigroup $h_0 = r_0 = \alpha_0=\gamma_0$ is isotopic
to $Z_k$ (otherwise after permuting the arguments we get case (2a)).
Denote $s*t \triangleq r_0(s,t)$ (recall that $s*0\equiv 0*s\equiv s$).
By the example $z=2$, $w=3$, we will prove that
$f(x,y,z,w)\equiv x+y+z+v$.
From (\ref{eq:W}) we have:
$$\varphi_{2,0}(x,y)\equiv x+(2*y).$$
From (\ref{eq:Z}), similarly to (\ref{eq:a0ab}), we have
$$\varphi_{2,3}(x,y)\equiv \varphi_{2,0}(\sigma (x),y) \equiv \sigma (x)+(2*y)$$
for some permutation $\sigma$. Fixing $y=0$, we get
$$f(x,0,2,3)\equiv \varphi_{2,3}(x,y) \equiv\sigma (x)+2;$$
but from (\ref{eq:Y}) we have
$$ f(x,0,2,3)\equiv (x+2)*3;
$$
so, $\sigma (x) \equiv (x+2)*3-2$ and
\begin{equation}\label{eq:23xy}
\varphi_{2,3}(x,y) \equiv ((x+2)*3)-2+(2*y).
\end{equation}
Similarly, from (\ref{eq:W}), (\ref{eq:Z}), and (\ref{eq:X}) we derive
$$\varphi_{2,3}(x,y)\equiv \varphi_{0,3}(x,\tau (y)) \equiv  (x*3)+ \tau (y) \equiv (x*3)-3+(2*(3+y)).$$
Equating these two representations of $\varphi_{2,3}(x,y)$,
 we get
$$((x+2)*3)-2+(2*y) \equiv (x*3)-3+(2*(3+y));$$
equivalently,
$$((x+2)*3)-(x*3) \equiv (2*(3+y))-(2*y)-1.$$
Since the right part does not depend on $x$, we have $((x+2)*3)-(x*3) \equiv \mathrm{const}$.
By Proposition~\ref{prozk21} we find $x*3 \equiv c\cdot x+d$ for some constants $c$ and $d$.
If $c\neq 1$, then the equation $c\cdot x+d = x$, i.e., $x*3=x*0$, has a solution, which contradicts the fact that $*$ is a quasigroup.
Hence, $c=1$. From $0*3 = 3$ we also find $d=3$.
So, $x*3 \equiv x+3$, and from (\ref{eq:23xy}) we see
that $f(x,y,2,3)\equiv x+y+2+3$.
Similarly $f(x,y,z,w)\equiv x+y+z+w$ for every $z$ and $w$.
\proofend

\section{Sublinear $n$-quasigroups}\label{s:Z5Z7}

We will say that an $n$-quasigroup is {\em sublinear}
if all its $2$-retracts are isotopic to $Z_k$.
There are two following facts, which currently have only
computational confirmation:

\begin{proposition}
All sublinear $4$-quasigroups of order $5$ are reducible.
All sublinear $3$-quasigroups of order $7$ are reducible.
\end{proposition}

The first fact was checked directly, using the classification
\cite{MK-W:small,MK-W:data}
of the $4$-quasigroups (latin $4$-cubes) of order $5$. All
the sublinear $3$-quasigroups $f$ of order $7$ were tested using the
following approach. We first may take without loss of generality
$f(0,y,z) \equiv y+z$ and $f(x,0,0) \equiv x$. Then, there are
$5!=120$ ways to choose the retract $f(x,y,0)$ isotopic to $Z_7$,
and $120$ ways to choose $f(x,0,z)$. Finally, we choose the
permutation $f(1,1,z)$, after which, by Corollary~\ref{cor:z5z5},
the sublinear $3$-quasigroup must be unequally reconstructible.

By induction, using Theorem~\ref{th}, we obtain the following.

\begin{corollary}\label{cor:subli}
All sublinear $n$-quasigroups of order $5$ are reducible provided $n\geq 4$; of order $7$, provided $n\geq 3$.
\end{corollary}

There is a unique, up to isotopy, irreducible sublinear $3$-quasigroup of order $5$.
It would be quite interesting to generalize Corollary~\ref{cor:subli} to other orders.
Taking into account Theorem~\ref{th}, if the order is prime,
it is sufficient to prove the reducibility for some fixed arity, say $3$ or $4$.



\section*{}
\section*{}

\providecommand\href[2]{#2} \providecommand\url[1]{\href{#1}{#1}}
 \def\DOI#1{{DOI}: \href{http://dx.doi.org/#1}{#1}}
\def\DOIURL#1#2{{DOI}: \href{http://dx.doi.org/#2}{#1}}
  \def\arXiv:#1 {{\small {arXiv}: \href{http://arXiv.org/abs/#1}{#1}}}

\end{document}


\begin{frontmatter}
\title{О связи между разделимостью $n$-арной
квазигруппы и разделимостью ее ретрактов\thanksref{kropot}
}
\author{{Д.~С.~Кротов}},
\ead{krotov@math.nsc.ru}
\author{{В.~Н.~Потапов}}
\ead{vpotapov@math.nsc.ru}
\address{ИМ СО РАН, пр. Акад. Коптюга 4, Новосибирск 630090, Россия}
\address{ММФ НГУ, ул. Пирогова 2, Новосибирск 630090, Россия}

\thanks[kropot]{Это авторский перевод статьи в {Discrete Math. 311(1) (2011) 58--66}, DOI:~\href{http://dx.doi.org/10.1016/j.disc.2010.09.023}{10.1016/j.disc.2010.09.023}.}

\begin{abstract}
 $n$-Арная операция $Q:\Sigma^n\to \Sigma$ называется {$n$-арной
квазигруппой} {порядка $|\Sigma|$}, если в уравнении
 $x_{0}=Q(x_1,\ldots , x_n)$
 значения любых $n$ переменных из $x_0$, \ldots ,
$x_n$ однозначно задают значение оставшейся переменной.
$n$-Арная квазигруппа $Q$ называется разделимой (перестановочно приводимой),
если $Q(x_1,\ldots,x_n)=P\left(R(x_{\sigma
(1)},\ldots,x_{\sigma (k)}), \linebreak[1] x_{\sigma
(k+1)},\ldots,x_{\sigma (n)}\right)$, где $P$ и $R$ суть
$(n{-}k{+}1)$-арная и $k$-арная квазигруппы, $\sigma$ --- перестановка и $1<k<n$.
$m$-Арная квазигруппа $R$ называется ретрактом $n$-арной квазигруппы $Q$,
если она может быть получена из $Q$ или одного из обращений $Q$ фиксацией
значений некоторых $n-m>0$ аргументов.

Доказано, что каждая неразделимая $n$-арная квазигруппа имеет неразделимый
 $(n{-}1)$- или $(n{-}2)$- арный ретракт,
 более того, если порядок конечный и простой, то она имеет неразделимый $(n{-}1)$-арный ретракт.
Этот факт используется чтобы показать, что все $n$-арные квазигруппы порядка $5$ или $7$,
все
бинарные ретракты которых изотопны группе $Z_5$ или $Z_7$, являются разделимыми при $n\geq 4$.
\begin{keyword}
$n$-арная квазигруппа\sep
ретракт\sep
разделимость\sep
приводимость\sep
латинский гиперкуб
\MSC
05B99\sep 20N15\sep 94B25
\end{keyword}
\end{abstract}
\end{frontmatter}

В данной работе мы доказываем признак разделимости $n$-арных квазигрупп
в терминах разделимости ретрактов
(разд.~\ref{s:main}--\ref{s:4-2}) и применяем этот признак к классификации
$n$"=арных квазигрупп порядков $5$ и $7$, все бинарные ретракты которых изотопны группам $Z_5$ и
$Z_7$, соответственно (разд.~\ref{s:Z5Z7}).

Признак следующий: если все $(n{-}1)$- и $(n{-}2)$-арные ретракты $n$"=арной квазигруппы $f$ разделимы,
то $f$ разделима; если все $(n{-}1)$-арные ретракты $n$"=арной квазигруппы $f$ конечного простого порядка разделимы,
то $f$ разделима. Чтобы установить этот факт, мы дополняем результат предыдущей работы \cite{RU:Kro:n-3}.
Следует заметить, что любой признак разделимости имеет смысл только при условии существования неразделимых квазигрупп.
Для различных порядков и арностей этот факт был доказан в работах
\cite{RU:BelSan.ru},\cite{RU:Frenkin.ru},\cite{RU:Borisenko.ru},\cite{RU:Glu76.ru},\cite{RU:Gold75.ru},\cite{RU:Gold76.ru},\cite{RU:AkGo},\cite{RU:KPS:ir}.
Единственность канонической декомпозиции разделимой $n$"=арной квазигруппы на группы и неразделимые квазигруппы меньшего порядка
доказана в \cite{RU:Cher.ru}.

\section{Определения}\label{s:1}

\begin{definition}
{
 $n$-Арная операция $Q:\Sigma^n\to \Sigma$,
где $\Sigma$ --- непустое множество,
называется \emph{$n$"=арной квазигруппой}, или \emph{$n$"=квазигруппой}
(\emph{порядка $|\Sigma|$}), если 
в уравнении $x_{0}=Q(x_1, \ldots , x_n)$ значения любых $n$  переменных из  $x_0$, $x_1$, \ldots , $x_n$ однозначно задают значение оставшейся переменной {\rm \cite{RU:Belousov.ru}}.
}
\end{definition}

Из соображений симметрии, для $n$"=квазигруппы $q$ мы также будем использовать предикатное обозначение
$q\langle\cdot\rangle$, определяемое следующим образом:
$$
q \langle x_0,x_1,\ldots,x_n\rangle \Leftrightarrow x_{0}=q(x_1,\ldots,x_n).
$$

\begin{definition}
{
 Если мы присвоим некоторые фиксированные значения $l\in \{1,\ldots,n\}$ переменным
 в предикате $Q\langle x_0, \ldots , x_n\rangle$, то полученный $(n-l+1)$-арный предикат соответствует некоторой $(n-l)$"=квазигруппе,
 которая называется \emph{ретрактом}, или \emph{$(n-l)$"=ретрактом} $n$"=квазигруппы $Q$.
Если при этом переменная $x_0$ остается незафиксированной, то ретракт называется \emph{главным}.
}
\end{definition}

\begin{definition}
{
  {\em Изотопией} будем называть набор из $n+1$ перестановок $\tau_i:\Sigma\to \Sigma$,
$i\in \{0,1,\ldots,n\} $.
$n$-Квазигруппы $f$ и $g$ называются {\em изотопными},
если для некоторой изотопии $\bar \tau=(\tau_0,\tau_1,\ldots ,\tau_n)$ имеем
$f(x_1,\dots,x_{n}) \equiv \tau_0^{-1} g(\tau_1 x_{1}, \ldots, \tau_{n} x_{n})$, т.\,е.
$f\langle x_0,x_1,\dots,x_{n}\rangle \equiv g\langle \tau_0 x_{0}, \tau_1 x_{1}, \ldots, \tau_{n} x_{n}\rangle$.
}
\end{definition}

\begin{definition}
{
 $n$-Квазигруппа $f$ называется {\em разделимой} ({\em перестановочно приводимой}),
если найдутся: целое число $m$ из промежутка $\{2,\ldots,n-1\}$,
$(n-m+1)$"=квазигруппа $h$, $m$"=квазигруппа $g$ и перестановка $\sigma: \{1,\ldots,n\} \to \{1,\ldots,n\}$ ---
такие, что
$$f(x_1,\ldots,x_{n}) \equiv h(g(x_{\sigma(1)},\ldots, x_{\sigma(m)}),
 x_{\sigma(m+1)},\ldots, x_{\sigma(n)}).$$
Если $n$"=квазигруппа не является разделимой, то она {\em неразделима}.
В частности, все бинарные
(как и $1$- и $0$-арные) квазигруппы неразделимы.
}
\end{definition}

\section{Основные результаты}\label{s:main}
Мы докажем (разд.~\ref{s:k2}--\ref{s:4-2}) следующие две леммы.
\begin{lemma}\label{l:k2}
Пусть все главные $3$- и $4$-ретракты $n$"=квазигруппы $f$ ($n\geq 5$) разделимы.
Тогда $f$ разделима.
\end{lemma}

Лемма~\ref{l:k2} была доказана в работе \cite{RU:KroPot:4} для случая порядка~$4$.
На самом деле, то доказательство можно рассматривать как упрощенную версию%
\footnote{Еще более слабое утверждение, доказательство которого тем не менее иллюстрирует ход рассуждений 
доказательства леммы~\ref{l:k2}, рассмотрено в \cite{RU:Kro:swi_sep.ru} для так называемой свитчинговой разделимости графов,
которая эквивалентна разделимости построенных по этим графам квазигруппам.}
доказательства леммы~\ref{l:k2}.
(В случае порядка $4$ все $2$"=квазигруппы изотопны коммутативным группам, таким образом, нет нужды
рассматривать неассоциативные $2$"=квазигруппы и некоммутативные группы, что существенно упрощает доказательство.)

\begin{lemma}\label{l:n-2}
Пусть $n\geq 4$.
Если $n$"=квазигруппа $f$ конечного простого порядка имеет неразделимый $(n{-}2)$"=ретракт
и все ее $(n{-}1)$"=ретракты разделимы, то $f$ разделима.
\end{lemma}

Эти леммы дополняют результат работы \cite{RU:Kro:n-3}, в которой утверждается, что если
максимальная арность $\kappa$ неразделимого ретракта $n$"=квазигруппы принадлежит промежутку $\{3,\ldots,n-3\}$,
то сама $n$"=квазигруппа является разделимой.
В самом деле, лемма~\ref{l:k2} решает случай $\kappa=2$, а лемма~\ref{l:n-2}
решает случай $\kappa=n-2$ для простых порядков. В итоге мы имеем следующее:

\begin{theorem}\label{th}
Пусть $f$ --- неразделимая $n$"=квазигруппа, $n\geq 4$.
Тогда $f$ имеет неразделимый $(n{-}1)$- или $(n{-}2)$-ретракт.
Более того, если при этом порядок $n$"=квазигруппы $f$ конечен и прост,
то $f$ имеет неразделимый $(n{-}1)$-ретракт.
\end{theorem}

Как показано в работе \cite{RU:Kro:n-2}, для каждой четной арности $n\geq 4$ и каждого порядка $q=4k$
 существует неразделимая $n$"=квазигруппа, у которой все $(n-1)$"=ретракты разделимы.
Таким образом, условия леммы~\ref{l:n-2} и последнего утверждения теоремы~\ref{th} не могут быть расширены на все порядки.
Случай нечетного $n$ и случай составного порядка $q \not\equiv 0 \bmod 4$
остаются неисследованными.

Хотя лемма~\ref{l:n-2} неприменима для порядка $4$, теорема~\ref{th}
(доказанная для этого порядка раньше) оказалась полезна при характеризации класса $n$"=арных квазигрупп порядка~$4$ \cite{RU:KroPot:4}.
Если возможна характеризация класса квазигрупп некоторого большего фиксированного порядка или некоторого подкласса, замкнутого относительно взятия ретракта,
то общая версия теремы~\ref{th} также может принести большую пользу, особенно для простого порядка:
рассуждая по индукции и доказывая некоторое предположение о структуре неразделимых $n$"=квазигрупп,
мы можем полагать, что $n$"=квазигруппа имеет неразделимый $(n-1)$"=ретракт, для которого доказываемое предположение верно.
Самый простой случай --- когда класс состоит только из разлелимых квазигрупп, примеры таких классов рассмотрены в разделе~\ref{s:Z5Z7}.

\section{Доказательство леммы~\ref{l:k2}}\label{s:k2}
Пусть $0$ --- некоторый элемент множества $\Sigma$.
Без потери общности мы можем рассмотреть $n$"=квазигруппу $f$, удовлетворяющую тождествам
$$f(0,\ldots,0,x_i,0,\ldots,0)\equiv x_i$$
для всех $i\in [n]$ и $x_i\in\Sigma$ (в противном случае,
квазигруппу можно привести к такому виду, применив подходящую
изотопию). Предположим, что все $3$-арные и $4$-арные ретракты
квазигруппы $f$ являются разделимыми.

Для обозначения бинарных ретрактов квазигруппы $f$ будем использовать символы операций
вроде $*$, $\star$, $\diamond$.
Если $\star$ --- бинарная операция, то через $\bar\star$ будем обозначать
сопряженную ей операцию, т.\,е. $x\star y\equiv y\,{\bar\star}\, x$.
Будем писать $f|i*(j\star k)$, если $i,j,k\in[n]$
попарно различны и
$$f(0,\ldots,0,x_i,0,\ldots,0,x_j,0,\ldots,0,x_k,0,\ldots,0)\equiv x_i*(x_j\star x_k)$$
(в общем случае индексы $i,j,k$ могут следовать в любом порядке),
где $*$ и $\star$ --- бинарные квазигруппы.
При этом если $*$ совпадает с $\star$ и является ассоциативной (т.\,е. группой),
то скобки будем иногда опускать.
Будем писать $f|i(jk)$, если $f|i*(j\star k)$ для некоторых $*$ и $\star$.
Будем писать $f|i[jk]$, если $f|i*(j\star k)$ для некоторых $*$ и $\star$
таких, что $\star\not\in\{*,\bar *\}$ или $*$ не ассоциативна.
Аналогичную запись будем использовать с двумя или четырьмя индексами.
Заметим следующее.
\begin{proposition} \label{p:1234}
Для любых попарно различных $i_1, i_2, i_3, i_4 \in \{1,\ldots,n\}$ выполнено либо
$f|i_a(i_b(i_c i_d))$, либо $f|(i_a i_b)(i_c i_d)$, где $\{a,b,c,d\}=\{1,2,3,4\}$.
\end{proposition}

\begin{proposition} \label{p:abc}
{\rm (a)} $f|i(jk)$  тогда и только тогда, когда $f|i(kj)$;\\
\hphantom{\rm (a)} $f|i[jk]$  тогда и только тогда, когда $f|i[kj]$.
\\
{\rm (b)} Для попарно различных $i,j,k\in [n]$ всегда выполнено хотя бы одно из
\\ \hphantom{\rm (a)}  $f|i(jk)$, $f|j(ik)$, $f|k(ij)$.
\\
{\rm (c)} $f|i[jk]$ тогда и только тогда,
\\ \hphantom{\rm (a)} когда верно $f|i(jk)$ и неверны ни $f|j(ik)$, ни $f|k(ij)$.
\\
{\rm (d)} Из $f|i(jk)$ и $f|j(ik)$ следует $f|i*k*j$ для некоторой ассоциативной $*$,
\\ \hphantom{\rm (a)} причем $*$ коммутативна тогда и только тогда, когда $f|k(ij)$.
\end{proposition}

\proofr
(a) вытекает из определений. Действительно, $f|i\star(j * k)$ эквивалентно
$f|i\star(k \,{\bar *}\, j)$.

(b) следует из разделимости $3$-ретрактов квазигруппы $f$.

(d) Предположим, что верно $f|i(jk)$ и $f|j(ik)$. Учитывая (a),
имеем $f|i*(k \star j)$ и $f|(i \circ k) \diamond j$, т.\,е.
\begin{equation}
\label{eq:assoc}
x * (y \star z) \equiv (x \circ y) \diamond z,
\end{equation}
для некоторых $*$, $\star$, $\circ$ и $\diamond$.
(Тождество (\ref{eq:assoc}) известно как обобщенное функциональное уравнение ассоциативности,
общее решение в квазигруппах см. \cite{RU:Belousov.ru}.)
Подставляя по очереди $x=0$, $y=0$ и $z=0$,
мы находим, что операции $*$, $\star$, $\circ$ и $\diamond$ --- одно и то же.
Тогда (\ref{eq:assoc}) означает ассоциативность $*$.
Очевидно, что из коммутативности будет следовать $f|k(ij)$.
Обратно, если $f|k{\scriptstyle\,\triangle\,}(i{\scriptstyle\,\square\,} j)$, то, подставляя ноль в
$x * y * z \equiv y{\scriptstyle\,\triangle\,}(x {\scriptstyle\,\square\,} z)$,
получаем $x * y \equiv y{\scriptstyle\,\triangle\,} x$
и $y * z \equiv y{\scriptstyle\,\triangle\,} z$, т.\,е. $*$ коммутативна.

(c) следует непосредственно из (a), (d) и определения обозначения $f|i[jk]$.
\proofend

\begin{proposition} \label{p:d}
{\rm (a)} Из $f|i[jk]$ и $f|j(kl)$ следует $f|i[j(kl)]$, в частности,
$i\neq l$, $f|i[jl]$, и $f|i(kl)$.
\\
{\rm (b)} Если $f|i*j*l$ и $f|j*k*l$, где $*$ ассоциативна и некоммутативна,
\\ \hphantom{\rm (a)} то $f|i*j*k*l$.
\end{proposition}
\proofr
(a) Пусть $f|i[jk]$  и $f|j(kl)$,
т.\,е. $f|i*(j \star k)$ и $f|j\circ (k \diamond l)$, где
$\star\not\in\{*,\bar *\}$ или $*$ неассоциативна.
Фиксируя нулями по очереди $i$-ю и $l$-ю переменную, мы видим, что
$\star$ и $\circ$ совпадают.
Легко проверить (с помощью предложения~\ref{p:1234}),
что единственный вариант для $4$-ретракта в переменных с номерами $i$, $j$, $k$, $l$ ---
$f|i*(j \star (k \diamond l))$.
Таким образом, имеем $f|i*(j \star l)$ и $f|i*(k \diamond l)$,
откуда $f|i[jl]$ и $f|i(kl)$ (заметим, что мы не можем гарантировать $f|i[kl]$, т.\,к.
$*$ может совпадать с $\diamond$ или $\bar \diamond$ и быть ассоциативной).

(b) Легко проверить по предложению~\ref{p:1234}, что любые другие варианты невозможны.
\proofend

\begin{definition}[вполне разделимая $n$"=квазигруппа]
{
$n$-Квазигруппа $f$ называется  {\em вполне разделимой},
если она разделимая и разделимы все ее главные ретракты арности больше $2$.
}
\end{definition}

\begin{lemma}\label{l:phif}
Предположим, что все главные $3$- и $4$-ретракты $n$-квазигруппы $f$ являются разделимыми.
Тогда найдется вполне разделимая $n$"=квазигруппа $\phi_f$, которая совпадает с $f$
на всех наборах с не более чем тремя ненулевыми элементами.
\end{lemma}

\proofr Будем рассуждать индукцией по $n$.
Мы покажем существование <<внутренней>> пары аргументов,
тогда <<стягивание>> этой пары позволит уменьшить $n$.

Будем говорить, что пара $\{a,b\}\subset \{1,\ldots,n\}$ является \emph{внутренней} (\emph{квази-внутренней}),
если для каждого $c\in \{1,\ldots,n\}\setminus\{a,b\}$ выполнено $f|c(ab)$
(соответственно, $f|c(ab)$ или $f|a*c*b$ для некоторой ассоциативной $*$).

(I). Сначала докажем существование квази-внутренней пары.

Рассмотрим последовательность $b,a_1, \ldots, a_l$ элементов из $[n]$, $1\leq l <n$,
удовлетворяющую следующему условию 
\begin{equation}\label{eq:baa}
f|a_i[a_{i+1} b] \mbox{ для всех $i=1,\ldots,l-1$}.
\end{equation}

(*). \emph{Мы утверждаем, что для всех $i$ и $j$, удовлетворяющих $1\leq i < j \leq l$,
выполнено $f|a_i(a_{j} b)$.} Действительно, для $j=i+1$ это так.
По предложению~\ref{p:d}(a), $f|a_i[a_{i+1} b]$ и $f|a_{i+1}(a_{j} b)$
влекут $f|a_i(a_{j} b)$. Таким образом, (*) доказывается индукцией по $j-i$.

В частности, (*) означает, что все $a_i$ попарно различны,
и длина последовательности ограничена числом $n$.
Рассмотрим максимальную по длине последовательность
 ${\mathbf a}:\,b,a_1, \ldots, a_L$, удовлетворяющую условию (\ref{eq:baa}) с $l=L$.

(**). \emph{Пара $b$, $a_L$ является квази"=внутренней.}
Действительно, учитывая предложение~\ref{p:abc} и определение
квази"=внутренней пары, достаточно исключить следующие два случая:

i) $f|a_L[bc]$ для некоторого $c$. Но это приводит к противоречию с максимальностью $L$,
поскольку последовательность $b,a_1, \ldots, a_L, a_{L+1}=c$ также удовлетворяет условию (\ref{eq:baa}) с $l=L+1$.

ii) $f|b[a_L c]$. В этом случае также можно показать, что  $L$ не максимально.
В самом деле, рассмотрим последовательность ${\mathbf a'}:\,c,a_1, \ldots, a_{L-1}, b, a_L$.
Имеем \underline{$f|b[a_L c]$}.
Если $L=1$, то доказывать нечего.
В противном случае, по предложению~\ref{p:d}(a), из $f|a_{L-1}[a_L b]$ и  $f|b[a_L c]$
следует
\underline{$f|a_{L-1}[b c]$}. 
Рассуждая индуктивно, из $f|a_i[a_{i+1} b]$ и $f|a_{i+1}(b c)$ имеем
\underline{$f|a_i[a_{i+1} c]$} и $f|a_i(b c)$ (последнее используется для очередного шага индукции)
для $i=L-1,L-2,\ldots,2,1$.
Таким образом, подчеркнутые выражения показывают, что последовательность
${\mathbf a'}$ удовлетворяет условию (\ref{eq:baa}) с $l=L+1$ и, следовательно, $L$ не максимально.
Полученное противоречие доказывает (**) и завершает (I).

(II). Теперь нашей целью является нахождение внутренней пары.
Воспользуемся последовательностью ${\mathbf a}$, построенной раннее.
По определению, единственная возможность, чтобы квази"=внутренняя пара
$b$, $a_L$ не являлась внутренней --- существование такого $d_1$,
что выполнено $f|a_L*d_1*b$, где $*$ --- ассоциативная и некоммутативная операция,
однозначно определяемая соотношением $f|a*b$. В этом случае

(***) \emph{мы утверждаем, что $f|a_{L-1}[d_1 b]$, т.\,е. мы можем заменить элемент $a_L$ элементом $d_1$ в последовательности ${\mathbf a}$.}
В самом деле, по предложению~\ref{p:d}(a), из $f|a_{L-1}[a_{L} b]$ и $f|a_L(d_1*b)$ вытекает
$f|a_{L-1}(d_1 b)$. Более того, операции в $f|a_{L-1}\star (d_1 * b)$ совпадают с операциями в $f|a_{L-1} \star (a_{L} b)$,
и (***) верно по определению.

Таким образом, мы получаем другую квази"=внутреннюю пару $b$, $d_1$.
В свою очередь,  либо она является внутренней, либо существует $d_2$ такой, что $f|d_1*d_2*b$.
Рассуждая таким образом, мы построим последовательность
${\mathbf d}:a_L=d_0,d_1,d_2,\ldots$ такую,
что $f|d_{i-1}*d_{i}*b$ и пара $b$, $d_i$ квази"=внутренняя.
Пользуясь предложением~\ref{p:d}(b), мы видим по индукции, что $f|d_j * d_i * b$ для любых $j<i$.
Это означает, что все $d_i$ попарно различны и последовательность $\mathbf d$ не может быть бесконечной,
т.\,е. на некотором шаге с номером $r$ мы получили внутреннюю пару $b$, $d_r$.

(III) Тривиально, утверждение леммы верно для $n = 4$. Предположим, что оно верно для $n=k-1$.
Рассмотрим случай $n=k$.
Как показано ранее, $f$ имеет внутреннюю пару аргументов. Без потери общности, это пара $\{n-1,n\}$.
Пусть $*$ --- соответствующая операция, т.\,е. $f|(n-1) * n$.
И пусть $(n-1)$"=квазигруппа $q$ получена из $f$ фиксацией нулем последнего аргумента.
Легко поверить непосредственно, что $n$"=квазигруппа $\phi_f$,
определяемая равенством
$\phi_f(x_1,\ldots,x_n)\triangleq \phi_q(x_1,\ldots,x_{n-2},(x_{n-1}*x_n))$,
удовлетворяет условию леммы.
\proofend

Последнее вспомогательное утверждение, которое нам понадобится, следующее:
\begin{proposition}[{\rm\cite[Theorem~1]{RU:KPS:ir}}] \label{p:24}
 Пусть $q$ и $g$ --- разделимые $4$"=квазигруппы, и пусть
\begin{eqnarray*}
  q(x,y,z,0)&\equiv& g(x,y,z,0), \\
  q(x,y,0,w)&\equiv& g(x,y,0,w), \\
  q(x,0,z,w)&\equiv& g(x,0,z,w), \\
  q(0,y,z,w)&\equiv& g(0,y,z,w).
\end{eqnarray*}
 Тогда $q\equiv g$.
\end{proposition}

\proofr[леммы~\ref{l:k2}]
По лемме~\ref{l:phif} существует вполне разделимая $n$"=квазигруппа $\phi_f$,
совпадающая с $f$ на $n$"=наборах, содержащих не более трех ненулевых элементов.
Пользуясь предложением~\ref{p:24}, легко доказать индукцией по числу ненулевых элементов в
$\bar x
$, что
$\phi_f(\bar x)=f(\bar x)$ for any $\bar x$.

Действительно, пусть $i,j,k,l$, $1\leq i < j<k<l\leq n$, ---
четыре номера, соответствующие ненулевым позициям набора $\bar x$.
Обозначим через $f_{\bar x; i,j,k,l}$ и $\phi_{\bar x; i,j,k,l}$ $4$-ретракты квазигрупп $f$ и $\phi_f$,
соответственно, у которых свободные аргументы имеют номера $i,j,k,l$, а остальные аргументы фиксированы
значениями из $\bar x$.
По предположению индукции значения $4$"=квазигрупп $f_{\bar x; i,j,k,l}$ и $\phi_{\bar x; i,j,k,l}$
совпадают на четверках, содержащих хотя бы один нуль, т.\,е.
удовлетворяют предположению предложения~\ref{p:24}. Поэтому они тождественно равны;
в частности,
$f_{\bar x; i,j,k,l}(x_i,x_j,x_k,x_l)=\phi_{\bar x; i,j,k,l}(x_i,x_j,x_k,x_l)$,
т.\,е. $f(\bar x)=\phi_f(\bar x)$.
\proofend

\section{Доказательство леммы~{\ref{l:n-2}}, случай $n>4$}\label{s:n-2}

The line of reasoning of the proof reminds the proof of Lemma~3 in \cite{RU:KroPot:4}.
But there are crucial differences that do not allow to unify these two statements:
in \cite[Lemma~3]{RU:KroPot:4} we deal with $n$-ary quasigroups of order $4$ and prove the
reducibility or semilinearity (some special property of $n$-ary quasigroups of order $4$);
while here we deal with a prime order and prove the reducibility.
To find a special property that would replace the semilinearity in a hypothetical general theorem
is an open problem.

Ход рассуждений в доказательстве текущего случая сильно напоминает
доказательство леммы~3 из \cite{RU:KroPot:4}. Но имеются и существенные
различия, которые не позволяют объединить эти два утверждения: в
\cite[Lemma~3]{RU:KroPot:4} мы имеем дело с $n$-квазигруппами порядка $4$
и доказываем разделимость или полулинейность (некоторое специальное
свойство $n$-квазигрупп именно порядка $4$), в то время как в
текущем разделе рассматривается только простой порядок и
доказывается разделимость. Специальное свойство, которое могло бы
обобщить полулинейность в гипотетической общей теореме подобного
рода, верной для любого порядка, --- открытая проблема.

\begin{proposition}[{\rm\cite[Proposition~5.1]{RU:KroPot:4}}] \label{p:canon}\label{RAQ:p:canon}
 Пусть разделимая $n$-квазигруппа $D$ ($n\geq 3$) имеет неразделимый $(n-1)$-ретракт
 $F\langle x_0,\ldots,x_{n-1}\rangle \equiv  D\langle x_0,\ldots,x_{n-1},0\rangle$.
Тогда найдутся  $i\in \{0,\ldots ,n\}$ и $2$-ква\-зи\-груп\-па $h$
такие, что $h(x,0)\equiv x$ и
\begin{equation}\label{RAQ:eq:canon}\label{eq:canon}
D\langle x_0,\ldots,x_{n}\rangle\equiv F\langle x_0,\ldots,
x_{i-1},h(x_i,x_{n}), x_{i+1},\ldots, x_{n-1}\rangle.
\end{equation}
\end{proposition}

\begin{proposition} \label{p:aut}
Пусть $f$ --- $n$-квазигруппа простого порядка, $n\geq 3$. Если
тождество
$$f\langle x_0,\ldots,x_n \rangle \equiv
f\langle
x_0,\ldots,x_{i-1},\mu(x_i),x_{i+1},\ldots,x_{j-1},\nu(x_j),x_{j+1},\ldots,x_n
\rangle$$ выполнено для некоторых различных $i,j\in {0,...,n}$ и
некоторой пары перестановок $(\mu,\nu)$, отличных от тождественной,
то $n$-квазигруппа $f$ разделима.
\end{proposition}
\proofr
 Без потери общности будем считать, что $i=1$, $j=2$. Положим
\begin{eqnarray*}
\alpha(x,y) &\triangleq& f(x,y,\bar 0), 
\\
\beta(x,\bar z) &\triangleq& f(x,0,\bar z),\qquad \bar z\triangleq(z_1,\ldots ,z_{n-2}),\\
\gamma(x) &\triangleq& f(x,0,\bar 0). 
\end{eqnarray*}
Заметим, что
$$\gamma^{-1}(\alpha(x,0))\equiv x.$$
Пусть пара перестановок $(\mu, \nu)$ удовлетворяет предположению
предложения. Тогда $\alpha(x,y) \equiv \alpha(\mu x,\nu y)$ и
$\alpha(x,y) \equiv \alpha(\mu^s x,\nu^s y)$ для всех натуральных
$s$. Из определения квазигруппы, следует, что равенства $x=\mu^s x$
и $y=\nu^s y$ верны или не верны одновременно (при тех же $s$).
Таким образом, легко видеть, что перестановки $\mu$ и $\nu$ состоят
из циклов одинаковой длины, а поскольку порядок простой, они
являются циклическими перестановками. Как следствие, верно
следующее:

(*) Для каждого $v\in\Sigma$ существуют перестановки
$\rho_v,\tau_v:\Sigma\to \Sigma$ такие, что $f(x,y,\bar z)\equiv
f(\rho_v x,\tau_v y,\bar z)$ и $\tau_v v=0$ (другими словами, группа
перестановок $\tau$ таких, что $f(x,y,\bar z)\equiv f(\rho x,\tau
y,\bar z)$ для некоторого $\rho$, действует транзитивно на $\Sigma$,
т.\,е. имеет только одну орбиту).

 Тогда
\begin{eqnarray*}
f(x,y,\bar z)
&\equiv&
f(\rho_y x,\tau_y y,\bar z)
\equiv
f(\rho_y x,0,\bar z)
\equiv
\beta(\rho_y x,\bar z)
\\&\equiv&
\beta(\gamma^{-1}\alpha(\rho_y x,0),\bar z)
\equiv
\beta(\gamma^{-1}\alpha(\rho_y x,\tau_y y),\bar z)
\equiv
\beta(\gamma^{-1}\alpha(x,y),\bar z)
\end{eqnarray*}
и, таким образом, $f$ является разделимой при условии $n\geq 3$.
\proofend

\begin{proof}[леммы~\ref{l:n-2}, случай $n>4$]
Пусть $C$ --- $n$-квазигруппа. Предположим, что все ее
$(n-1)$-ретракты разделимы, а некоторый $(n-2)$-ретракт $E$
неразделим. Без потери общности будем считать, что
$$ E\langle x_0,\ldots,x_{n-2} \rangle \equiv C \langle x_0,\ldots,x_{n-2},0,0 \rangle. $$
Мы будем использовать следующие обозначения для ретрактов
$n$"=квазигруппы $C$:

\begin{eqnarray}
  E_{a,b}\langle x_0,\ldots,x_{n-2} \rangle & \triangleq & C \langle x_0,\ldots,x_{n-2},a,b \rangle, \label{eq:E}\\
  A_b\langle x_0,\ldots,x_{n-2},y \rangle & \triangleq & C\langle x_0,\ldots,x_{n-2},y,b \rangle, \label{eq:A}\\
  B_a\langle x_0,\ldots,x_{n-2},z \rangle & \triangleq & C\langle x_0,\ldots,x_{n-2},a,z \rangle. \label{eq:B}
\end{eqnarray}

(*) Поскольку квазигруппа $A_0$ разделима и фиксация $y:=0$ приводит
к неразделимому ретракту $E$, по предложению~\ref{p:canon} мы имеем
\begin{equation}\label{eq:A0}
A_0\langle x_0,\ldots,x_{n-2},y \rangle
\equiv
E\langle x_0,\ldots,x_{i-1},h(x_i,y),x_{i+1},\ldots,x_{n-2} \rangle
\end{equation}
для некоторых $i\in \{0,\ldots n-2\} $ и $2$-квазигруппы $h$ такой,
что $h(x_i,0)\equiv x_i$.

Из тождества (\ref{eq:A0}) видно, что все ретракты $E_{a,0}$,
$a\in\Sigma$, изотопны $E$. Аналогично получаем, что

(**) все ретракты $E_{a,b}$, $a,b\in\Sigma$, изотопны $E$.

Тогда, аналогично рассуждениям (*), для каждого $b\in\Sigma$
выполнено
$$A_b \langle x_0,\ldots,x_{n-2},y \rangle
\equiv
E_{0,b} \langle x_0,\ldots,x_{i_b-1},h_b(x_{i_b},y),x_{i_b+1},\ldots,x_{n-2} \rangle
$$
для некоторых $i_b\in \{0,\ldots n-2\} $ и $2$-квазигруппы $h_b$
такой, что $h_b(x,0)\equiv x$.

(***) Мы утверждаем, что номер $i_b$ не зависит $b$. Действительно,
предположим, например, что $i_1=0$ и $i_2=1$, т.\,е.
\begin{eqnarray*}
A_1 \langle x_0,\ldots,x_{n-2},y \rangle
&\equiv&
E_{0,1} \langle h_1(x_0,y),x_1,x_2,\ldots,x_{n-2} \rangle,\\
A_2 \langle x_0,\ldots,x_{n-2},y \rangle
&\equiv&
E_{0,2} \langle x_0,h_2(x_1,y),x_2,\ldots,x_{n-2} \rangle.
\end{eqnarray*}
Тогда, фиксируя $x_0$ в первом случае, мы получаем ретракт,
изотопный $E$; фиксируя $x_0$ во втором случае, мы получаем
разделимый ретракт (напомним, что $n\geq 5$). Но, аналогично
рассуждениям (**), эти два ретракта изотопны, что противоречит
разделимости $E$ и доказывает (***).

Без потери общности положим $i_b=0$, т.\,е.
\begin{equation}\label{eq:A0E}
A_b \langle x_0,x_1,\tilde x_2,y \rangle \equiv E_{0,b} \langle
h_b(x_0,y),x_1,\tilde x_2 \rangle;
\end{equation}
{здесь и далее
$\tilde x_2\triangleq (x_2,\ldots,x_{n-2})$}. Аналогично, мы можем без потери общности считать, что либо
\begin{equation}\label{eq:B0E}
B_a \langle x_0,x_1,\tilde x_2,z \rangle \equiv E_{a,0} \langle
g_a(x_0,z),x_1,\tilde x_2 \rangle,
\end{equation}
либо
\begin{equation}\label{eq:B1E}
B_a \langle x_0,x_1,\tilde x_2,z \rangle \equiv E_{a,0} \langle
x_0,g_a(x_1,z),\tilde x_2 \rangle,
\end{equation}
где $2$-квазигруппы $g_a$ удовлетворяет условию $g_a(x,0)\equiv x$.

Предполагая, что выполнено тождество (\ref{eq:B0E}), мы выводим
\begin{eqnarray}
C \langle x_0,x_1,\tilde x_2,y,z\rangle
&\stackrel{\mathrm{(\ref{eq:A})}}{\equiv}& \nonumber
A_z\langle x_0,x_1,\tilde x_2,y \rangle
\\
  &\stackrel{\mathrm{(\ref{eq:A0E})}}{\equiv}& \nonumber
  E_{0,z}\langle h_z(x_0,y),x_1,\tilde x_2\rangle
\\
  &\stackrel{\mathrm{(\ref{eq:E}),(\ref{eq:B})}}{\equiv}& \label{eq:4-5}
  B_{0}\langle h_z(x_0,y),x_1,\tilde x_2, z\rangle
\\
  &\stackrel{\mathrm{(\ref{eq:B0E})}}{\equiv}& \nonumber
  E_{0,0} \langle g_0(h_z(x_0,y),z),x_1,\tilde x_2 \rangle,
\end{eqnarray}
что означает разделимость $C$, поскольку $f(x,y,z) \triangleq
g_0(h_z(x,y),z)$ обязана быть $3$"=квазигруппой (фиксируя $x_1$ и
$\tilde x_2$, мы видим что $f$ является ретрактом $C$). Таким
образом, остается рассмотреть случай (\ref{eq:B1E}). Рассмотрим два подслучая:

Случай 1: $g_a$ не зависит от $a$; обозначим $g \triangleq g_0=g_a$
для всех $a\in\Sigma$. Тогда, повторяя три шага из (\ref{eq:4-5}),
мы выводим
\begin{eqnarray*}
C\langle x_0,x_1,\tilde x_2,y,z \rangle
  &\stackrel{\mathrm{(\ref{eq:B})}}{\equiv}&
B_y\langle x_0,x_1,\tilde x_2,z \rangle
  \\&\stackrel{\mathrm{(\ref{eq:B1E})}}{\equiv}&
E_{y,0} \langle x_0,g(x_1,z),\tilde x_2,z \rangle
  \\&\stackrel{\mathrm{(\ref{eq:E}),(\ref{eq:A})}}{\equiv}&
A_0 \langle x_0,g(x_1,z),\tilde x_2,y \rangle
\end{eqnarray*}
и видим, что $C$ разделимая.

Случай 2: для некоторого фиксированного $a$ имеем $g_0\neq g_a$;
обозначим $s_i(x)\triangleq g_0(x,i)$, $t_i(x)\triangleq g_a(x,i)$ и
$r_i(x) \triangleq h_i(x,a)$. Из (\ref{eq:A0E}) мы видим, что
\begin{eqnarray}
E_{a,0}\langle x_0,x_1,\tilde x_2 \rangle
&\equiv & \label{eq:E-1}
E_{0,0}\langle r_0(x_0),x_1,\tilde x_2 \rangle,
\\
E_{a,b}\langle x_0,x_1,\tilde x_2 \rangle
&\equiv & \label{eq:E-2}
E_{0,b}\langle r_b(x_0),x_1,\tilde x_2 \rangle.
\end{eqnarray}
Из (\ref{eq:B1E}) мы видим, что
\begin{eqnarray}
E_{0,b}\langle x_0,x_1,\tilde x_2 \rangle
&\equiv & \label{eq:E-3}
E_{0,0}\langle x_0,s_b(x_1),\tilde x_2 \rangle,
\\
E_{a,b}\langle x_0,x_1,\tilde x_2 \rangle
&\equiv & \label{eq:E-4}
E_{a,0}\langle x_0,t_b(x_1),\tilde x_2 \rangle.
\end{eqnarray}
Применяя последовательно (\ref{eq:E-3}), (\ref{eq:E-2}),
(\ref{eq:E-4}) и (\ref{eq:E-1}), мы находим, что для каждого $b$
ретракт $E=E_{0,0}$ удовлетворяет
$$
E\langle x_0,x_1,\tilde x_2 \rangle
\equiv E_{0,b}(\ldots)
\equiv E_{a,b}(\ldots)
\equiv E_{a,0}(\ldots)
\equiv E\langle r_b^{-1} r_0 x_0,s_b^{-1} t_b x_1,\tilde x_2 \rangle.
$$
По предложению~\ref{p:aut}, неразделимость $E$ означает, что
$s_b^{-1} t_b=Id$ для каждого $b$. Но это противоречит предположению
$g_0\neq g_a$ и, следовательно, случай 2 недопустим.
\end{proof}

\section{Доказательство леммы~\ref{l:n-2}, случай $n=4$}\label{s:4-2}

Пусть $\Sigma=\{0,1,\ldots,k-1\}$. В настоящем разделе
важную роль играет $2$"=квазигруппа $x+y\bmod k$,
которую будем обозначать через $Z_k$.
Все арифметические операции (сложение, умножение)
будем производить по модулю $k$, где $k$  --- порядок рассматриваемых квазигрупп.

\begin{proposition}\label{prozk1}
   Пусть $f$ --- $2$-квазигруппа простого порядка $k$.
   Если найдутся нетождественные перестановки
   $\nu, \mu:\Sigma\to\Sigma$ такие, что либо $f(\mu (x),\nu(y))\equiv f(x,y)$, либо
$f(\mu (x), y)\equiv \nu (f(x,y))$, либо
$f( x,\nu (y))\equiv \mu (f(x,y))$,
то квазигруппа $f$ изотопна $Z_k$.
\end{proposition}
\proofr
Рассмотрим только случай $f(\mu (x),\nu (y))\equiv f(x,y)$,
поскольку остальные случаи сводятся к нему рассмотрением квазигрупп $f^{\langle 1 \rangle}$ и $f^{\langle 2 \rangle}$,
полученных из $f$ обращением по первому и второму аргументу.

Сначала заметим, что перестановки $\mu$ и $\nu$ обязаны быть циклами на $k$
элементах (аналогичное утверждение содержится в доказательстве предложения~\ref{p:aut}).

Итак, докажем, что $f$ изотопна $Z_k$.
Найдутся перестановки $\alpha$ и $\beta$ такие, что
$\alpha\mu\alpha^{-1} = \epsilon$ и $\beta\nu\beta^{-1} = \epsilon^{-1}$,
где $\epsilon$ соответствует добавлению единицы $\bmod\ k$
(чтобы быть точным, перестановку $\alpha$ можно рекурсивно определить правилами
$\alpha(0)\triangleq 0$, $\alpha(\mu (x))\triangleq \alpha (x) +1$; аналогично $\beta$).
Тогда квазигруппа $g(x,y) \equiv f(\alpha^{-1} (x), \beta^{-1} (y))$ удовлетворяет
$g(\epsilon (x),\epsilon^{-1}(y)) \equiv g(x,y)$. Далее, определим
$h(x,y) \equiv \gamma^{-1}(g( x,  y))$,
где $\gamma(x) \triangleq g(x,0)$.
Теперь убедимся, что мы получили $Z_k$:
$
h(x,y)\equiv h(\epsilon^y (x), \epsilon^{-y} (y))\equiv h(x+y,0)\equiv x+y.
$
 \proofend

\begin{proposition}\label{prozk21}
   Пусть $|\Sigma|$ --- простое число.
   Предположим, что $\tau (x + a)\equiv \tau (x) +b$
    для некоторых $a,b\in \Sigma$ и перестановки $\tau:\Sigma\to\Sigma$.
    Тогда либо $\tau (x) \equiv c \cdot x +d$
   для некоторых  $c,d\in \Sigma$, либо $a=b=0$.  \end{proposition}

\proofr Если $a=0$, то очевидно $b=0$. Пусть $a\neq 0$.
Из $\tau (x + a)\equiv \tau (x) +b$ мы получаем по индукции, что
$\tau (y +
i\cdot a)\equiv \tau(y)+i\cdot b$ для каждого $i\in\Sigma$.
Далее,
$$\tau (x)\equiv \tau(0+(x/a)\cdot a) \equiv  \tau (0) + (x/a)\cdot b \equiv (b/a)\cdot x + \tau (0).$$
 \proofend

\begin{proposition}\label{prozk3} Пусть $f$ --- квазигруппа, изотопная $Z_k$.
Если $f(0,x)\equiv f(x,0)\equiv x$ и $f(x,1)\equiv x+1$, то $f(x,y)\equiv x+y$.
\end{proposition}

\proofr Имеем
$\alpha (f(x,y))\equiv
\beta (x) +\gamma (y)$
для перестановок $\alpha,\beta,\gamma:\Sigma\to\Sigma$.
Из $f(0,y)\equiv y$  получаем
\begin{equation} \label{eq:g=a-b}
\gamma(y) \equiv \alpha(y) - \beta( 0).
\end{equation}
Из $f(x,0)\equiv x$ получаем
\begin{equation} \label{eq:b=a-g}
\beta (x) \equiv \alpha (x) - \gamma (0).
\end{equation}
Из $f(x,1)\equiv x+1$ получаем
$\alpha(x+1) \equiv \beta (x)+\gamma (1) \equiv
\alpha (x) - \gamma (0)+\gamma (1)$.
Тогда по индукции
$\alpha (z)\equiv z\cdot (\gamma(1)-\gamma(0))+\alpha(0)$.
Следовательно
\begin{equation} \label{eq:a=a+a-a}
\alpha (x+y)\equiv\alpha (x)+\alpha (y)-\alpha(0).
\end{equation}
Пользуясь тождествами (\ref{eq:g=a-b})--(\ref{eq:a=a+a-a}), выводим
$$\alpha (f(x,y)) \equiv
\beta(x)+\gamma(y) \equiv
\alpha (x) - \gamma (0)+ \alpha(y) - \beta( 0) \equiv
\alpha (x) +\alpha(y) -\alpha( 0) \equiv \alpha(x+y),
$$
что доказывает тождество $f(x,y)\equiv x+y$.
 \proofend

\begin{corollary} \label{cor:z5z5} Пусть $k$ --- простое целое число,
пусть $2$"=квазигруппы $f$ и $g$ изотопны $Z_k$,
и пусть для некоторых фиксированных $a$, $a'$, $b$, $b'$, $a\neq b$,
выполнены тождества
$f(x,a)\equiv g(x,a')$ и $f(x,b)\equiv g(x,b')$.
Тогда $f(x,y)\equiv g(x,\pi(y))$ для некоторой перестановки $\pi$.
\end{corollary}
\proofr
Предложение~\ref{prozk3} является частным случаем рассматриваемого утверждения.
Мы сведем общий случай к этому частному.
Заметим, что утверждение эквивалентно аналогичному утверждению для
квазигрупп $f'$ и $g'$, полученных из $f$ и $g$
применением некоторой общей изотопии.
С учетом этого наблюдения, мы можем сделать следующие предположения:

1) без потери общности считать, что $g(x,y)\equiv x+y$;

2) пользуясь изотопией вида $(\mathrm{Id},\alpha,\alpha^{-1})$, где
$\alpha$ --- добавление константы, мы можем привести коэффициент $a'$ к нулю; далее считаем что $a'=0$;

3) пользуясь изотопией вида  $(\beta,\beta,\beta)$,
где $\beta$ --- умножение на константу, мы можем привести коэффициент $b'$ к единице, т.\,е. можно считать, что  $b'=1$
(в этом месте мы используем простоту числа $k$).

Более того, мы можем заменить $f$ квазигруппой $f''$, где $f''(x,y) \equiv f(x,\gamma(y))$ для некоторой перестановки $\gamma$.
Эту перестановку можно выбрать таким образом, что $f''(0,y) \equiv y$
(поскольку $a'=0$ и $b'=1$, мы также получим $a=0$, $b=1$).
Тогда утверждение следует из предложения~\ref{prozk3}.
\proofend

Следующее простое предложение мы оставляем без доказательства,
аналогичное утверждение (предложение~\cite[Proposition~2]{RU:Kro:n-3}) уже было доказано ранее.

\begin{proposition}\label{proraz4} Пусть $f$ --- $4$-квазигруппа конечного порядка,
и пусть все ее $3$"=ретракты разделимы. Тогда все ее $2$"=ретракты, полученные фиксацией
одних и тех же аргументов разными значениями, попарно изотопны.\end{proposition}

\begin{proposition}\label{proraz6} Пусть $f$ --- $4$-квазигруппа конечного порядка $k$,
и пусть существуют $2$"=квазигруппы $q$, $r$, $g_z$, $h_z$, $z\in\Sigma$,
такие, что
$$ f(x_1,x_2,x_3,z)\equiv g_z(x_3,h_z(x_1,x_2))\qquad \mbox{и}\qquad f(0,x_2,x_3,z)\equiv q(x_3,r(z,x_2)).$$
Тогда $4$"=квазигруппа $f$ разделима. \end{proposition}

\proofr
Подставляя $x_1=0$, получаем
$g_z(x_3,h_z(0,x_2))\equiv q(x_3,r(z,x_2))$.
Следовательно, для каждого $z$ имеем тождество $g_z(x_3,x_2)\equiv q(x_3,\pi_z(x_2))$ с некоторой перестановкой $\pi_z$.
Тогда $f(x_1,x_2,x_3,z)\equiv q(x_3,\pi_z^{-1}(h_z(x_1,x_2)))$, откуда следует разделимость
$f$.
\proofend

\proofr[леммы~\ref{l:n-2}, случай $n=4$]
Заметим, что любая $4$-квазигруппа порядка $k=2$ или $k=3$ изотопна итерированной группе $Z_k$, и поэтому разделима.
Пусть $k>3$.
Пусть $f:\Sigma^4\to\Sigma$ --- $4$"=квазигруппа и пусть все ее
$3$"=ретракты разделимы.
Для каждого $i=1,2,3,4$ и каждого $a\in \Sigma$ определим
\def\uot{\mathrm{out}}
$\uot(i,b)$ как порядковый номер аргумента $4$"=квазигруппы  $f$, который является  \emph{внешним}
в разложении ретракта квазигруппы $f(x_1,x_2,x_3,x_4)$, соответствующего фиксации переменной $x_i$ значением $b$.
Например, $\uot(1,2)=3$ означает $f(2,x_2,x_3,x_4)\equiv g(x_3,h(x_2,x_4))$ для некоторых $g$ и $h$.
Если соответствующий ретракт допускает два или три разложения с различными внешними аргументами,
то мы для определенности возьмем в качестве $\uot(i,b)$
минимальное из допустимых значений.

По мощностным соображениям, легко убедиться, что всегда справедливо по крайней мере одно из следующих двух утверждений:

(1) Имеются две различные позиции $i,j\in\{1,2,3,4\}$ и значения $a,c,d\in\Sigma$, $c\neq d$, такие, что
$\uot(i,a)=\uot(j,c)=\uot(j,d)$.

(2)  Функция $\uot(i) \triangleq \uot(i,b)$ не зависит от $b$ и взаимнооднозначна по отношению к $i$.
Более того, мы можем считать, что $\uot(i)$ --- циклическая перестановка,
 поскольку в противном случае для обращения $4$"=квазигруппы $f$ в любой позиции (скажем, в первой) выполнено утверждение (1).

Пусть верно (1).
Рассмотрим, без потери общности, следующие разложения ($j=3$, $i=4$, $a=c=0$, $d=1$):
\begin{eqnarray}
f(x,y,0,w)&\equiv& g_0(y,h_0(x,w)),\label{eq:f0}\\
f(x,y,1,w)&\equiv& g_1(y,h_1(x,w)),\label{eq:f1}\\
f(x,y,z,0)&\equiv& g(y,h(x,z))\label{eq:fx}.
\end{eqnarray}
Рассмотрим ретракты $\varphi_{z,w}(x,y)\triangleq
f(x,y,z,w)$. Дальнейшие рассуждения разделим на два подслучая.

(1a) \emph{Квазигруппа $\varphi_{0,0}$ изотопна $Z_k$}.

Подставляя $z=0$ и $w=0$ в тождествах (\ref{eq:f0}) и (\ref{eq:fx}),
мы видим, что квазигруппы $g_0$ и $g$ изотопны $Z_k$;
более того, $g_0(y,t)\equiv g(y,\tau(t))$ для некоторой перестановки $\tau$.
То же самое верно и для $g_1$. Таким образом, мы можем без потери общности считать, что
(\ref{eq:f0}) и (\ref{eq:f1}) выполняются с $g_0$ и $g_1$, равными $g$:
\begin{equation}\label{eq:g0g1g}
g_0=g_1=g.
\end{equation}

Теперь рассмотрим $2$-ретракты $q_{x,w}(y,z)\equiv f(x,y,z,w)$.
Из (\ref{eq:fx}) видно: если $w=0$, то квазигруппа $q_{x,w}$ изотопна $g$,
и следовательно изотопна $Z_k$. По предложению~\ref{proraz4}, это верно для всех $x$ и $w$.
Далее, из (\ref{eq:f0}) и (\ref{eq:g0g1g}) имеем: если $x$ и $w$ фиксированы,
то
$q_{x,w}(y,0)\equiv g(y,H_{0,x,w})$ для некоторой константы $H_{0,x,w}$ (точнее, $H_{0,x,w}\triangleq h_0(x,w)$).
Аналогично,
$q_{x,w}(y,1)\equiv g(y,H_{1,x,w})$ (где $H_{1,x,w}\triangleq h_1(x,w)$).
По следствию~\ref{cor:z5z5} имеем $q_{x,w}(y,z)\equiv g(y,\pi_{x,w}(z))$
для некоторой перестановки $\pi_{x,w}$. Но это означает, что  $f(x,y,z,w)\equiv g(y,\pi(x,z,w))$,
где
$\pi(x,z,w) \triangleq \pi_{x,w}(x)$, что доказывает разделимость $4$"=квазигруппы $f$.

(1b) \emph{Квазигруппа $\varphi_{0,0}$ не изотопна $Z_k$}. Покажем
следующее:

(*) \emph{все ретракты $4$"=квазигруппы $f$, полученные фиксацией четвертого аргумента,
имеют общий (второй) внешний аргумент}. От противного рассмотрим, для примера,
альтернативное разложение
\begin{equation}
f(x,y,z,2)\equiv g'(x,h'(y,z))\label{eq:fx'}.
\end{equation}

Из (\ref{eq:f0}) получаем
$\varphi_{0,2}(x,y) \equiv \varphi_{0,0}(\pi_0^{-1}(\pi_2(x)),y)$,
где $\pi_i(x) \equiv h_0(x,i)$ (эта явная формула для $\pi_i$ нам не понадобится).
Аналогично, из (\ref{eq:fx'}) получаем
$\varphi_{1,2}(x,y) \equiv \varphi_{0,2}(x,\nu_0^{-1}(\nu_1(y)))$,
где $\nu_i(y) \triangleq h'(y,i)$ (заметим, что $\nu_0\neq\nu_1$).
Аналогично, мы выразим $\varphi_{1,0}$ через $\varphi_{1,2}$, пользуясь (\ref{eq:f1}).
И, наконец, выразим $\varphi_{0,0}$ через $\varphi_{1,0}$, пользуясь (\ref{eq:fx}).
В результате получим
$\varphi_{0,0}(x,y)=\varphi_{0,0}(\mu (x),\nu (y))$ для некоторой перестановки $\mu$ и
 $\nu = \nu_0^{-1} \nu_1 \neq\mathrm{Id}$.
Пользуясь предложением~\ref{prozk1} мы заключаем, что $\varphi_{0,0}$ изотопна $Z_k$,
что противоречит условию рассматриваемого подслучая.
Аналогично, другой способ разложения $f(x,y,z,2)$ (с внешней $x$)
также приводит к противоречию. Более того, $2$ можно заменить на любое ненулевое значение.
Утверждение (*) доказано.

Теперь разделимость $4$"=квазигруппы $f$ следует непосредственно из предложения~\ref{proraz6}.

Пусть выполнено (2). Без потери общности пусть $\uot(3)=2$,
$\uot(2)=4$, $\uot(4)=1$, $\uot(1)=3$. Т.\,е.
\begin{eqnarray}
f(x,y,z,w)
 &\equiv& g_z(h_z(x,w),y)  \label{eq:Z} \\
 &\equiv& q_w(x,r_w(z,y))  \label{eq:W} \\
 &\equiv& \alpha_x(z,\beta_x(w,y))  \label{eq:X} \\
 &\equiv& \gamma_y(\delta_y(x,z),w). \label{eq:Y}
\end{eqnarray}
для некоторых квазигрупп $g_t$, $h_t$, $q_t$, $r_t$, $\alpha_t$,
$\beta_t$, $\gamma_t$, $\delta_t$, $t\in \Sigma$.
Без потери общности мы можем полагать, что
$f(0,0,0,v) \equiv f(0,0,v,0) \equiv  f(0,v,0,0) \equiv f(v,0,0,0) \equiv v$.
Более того, квазигруппы
$g_0$, $h_0$, $q_0$, $r_0$, $\alpha_0$, $\beta_0$,
$\gamma_0$, $\delta_0$
могут быть выбраны таким образом, что для любой квазигруппы $p$ из них верно
$p(0,v)\equiv p(v,0) \equiv v$ (см., например, \cite[предложение~12]{RU:PotKro:asymp.ru}).
Фиксируя в (\ref{eq:Z})--(\ref{eq:Y}) две переменные нулями, мы выводим
$g_0 = q_0 = \beta_0=\delta_0$ и
$h_0 = r_0 = \alpha_0=\gamma_0$. Дальнейшие рассуждения различные
в зависимости от того, изотопны все эти квазигруппы $Z_k$ или нет.

Мы снова рассмотрим ретракты
$\varphi_{z,w}(x,y)\triangleq f(x,y,z,w)$ и разобьем рассуждения на два подслучая.

(2a)
{\em Квазигруппы $\varphi_{z,w}$ (как и $g_0$, $q_0$, $\beta_0$, $\delta_0$) не изотопны $Z_k$.}
Рассмотрим ретракты $\varphi_{0,0}$, $\varphi_{0,b}$,
$\varphi_{a,0}$, $\varphi_{a,b}$ для некоторых  $a$ и $b$ из $\Sigma\setminus\{0\}$.
Из (\ref{eq:Z}) имеем
\begin{eqnarray*}
\varphi_{0,0}(x,y)&\equiv& f(x,y,0,0) \equiv
g_0(h_0(x,0),y),
\\
\varphi_{0,b}(x,y)&\equiv& f(x,y,0,b) \equiv
g_0(h_0(x,b),y).
\end{eqnarray*}
Таким образом,
\begin{equation}\label{eq:000b}
\varphi_{0,b}(x,y) \equiv \varphi_{0,0}(\sigma_{0,b}(x),y)
\end{equation}
для некоторой перестановки $\sigma_{0,b}$ (точнее,
$\sigma_{0,b}(x) \triangleq h_0(h_0^{-1}(x,0),b)$\,).
Аналогично, из (\ref{eq:Z}) и (\ref{eq:W}) мы выводим
\begin{eqnarray}
\label{eq:a0ab}
\varphi_{a,b}(x,y) \equiv \varphi_{a,0}(\sigma_{a,b}(x),y) \\
\label{eq:00a0}
\varphi_{a,0}(x,y) \equiv \varphi_{0,0}(x,\tau_{a,0}(y)) \\
\label{eq:0bab}
\varphi_{a,b}(x,y) \equiv \varphi_{0,b}(x,\tau_{a,b}(y))
\end{eqnarray}
для некоторых $\sigma_{a,b}$, $\tau_{a,0}$, $\tau_{a,b}$.
Мы покажем, что $\sigma_{0,b}=\sigma_{a,b}$
и $\tau_{a,0}=\tau_{a,b}$. Действительно, из тождеств (\ref{eq:000b})--(\ref{eq:0bab}) имеем
$$
\varphi_{0,0}(x,y)\equiv
\varphi_{0,0}(\sigma_{0,b}^{-1}(\sigma_{a,b}(x)), \tau_{a,b}(\tau^{-1}_{a,0}(y))),
$$
и из предложения~\ref{prozk1},
учитывая условие подслучая (2a),
мы заключаем, что
$\sigma_{0,b}^{-1}\sigma_{a,b}
=\tau_{a,b}\tau^{-1}_{a,0}=\mathrm{Id}$.

Теперь имеем
$$
f(x,y,z,w) \equiv \varphi_{z,w}(x,y)
\equiv \varphi_{0,0}(\sigma_{0,w} (x),\tau_{z,0} (y)),
$$
что означает разделимость $4$"=квазигруппы $f$.

(2b) {\em Квазигруппы $\varphi_{z,w}$ изотопны $Z_k$}.
Без потери общности предположим, что квазигруппы $g_0$, $q_0$, $\beta_0$ и $\delta_0$ совпадают с $Z_k$.
Мы также можем считать, что квазигруппа $h_0 = r_0 = \alpha_0=\gamma_0$ изотопна $Z_k$
(в противном случае переставив аргументы мы получаем подслучай (2a)).
Обозначим $s*t \triangleq r_0(s,t)$ (напомним, что $s*0\equiv 0*s\equiv s$).
На примере $z=2$, $w=3$ мы покажем, что
$f(x,y,z,w)\equiv x+y+z+v$.
Из (\ref{eq:W}) имеем:
$$\varphi_{2,0}(x,y)\equiv x+(2*y).$$
Из (\ref{eq:Z}), аналогично (\ref{eq:a0ab}), имеем
$$\varphi_{2,3}(x,y)\equiv \varphi_{2,0}(\sigma (x),y) \equiv \sigma (x)+(2*y)$$
для некоторой перестановки $\sigma$. Фиксируя $y=0$, получаем
$$f(x,0,2,3)\equiv \varphi_{2,3}(x,y) \equiv\sigma (x)+2.$$
Но из (\ref{eq:Y}) мы имеем
$$ f(x,0,2,3)\equiv (x+2)*3.
$$
Таким образом, $\sigma (x) \equiv (x+2)*3-2$, и
\begin{equation}\label{eq:23xy}
\varphi_{2,3}(x,y) \equiv ((x+2)*3)-2+(2*y).
\end{equation}
Аналогично, из (\ref{eq:W}), (\ref{eq:Z}), и (\ref{eq:X}) мы выводим
$$\varphi_{2,3}(x,y)\equiv \varphi_{0,3}(x,\tau (y)) \equiv  (x*3)+ \tau (y) \equiv (x*3)-3+(2*(3+y)).$$
Приравнивая эти два представления квазигруппы $\varphi_{2,3}(x,y)$,
 получаем
$$((x+2)*3)-2+(2*y) \equiv (x*3)-3+(2*(3+y)),$$
эквивалентно,
$$((x+2)*3)-(x*3) \equiv (2*(3+y))-(2*y)-1.$$
Поскольку правая часть не зависит от $x$, имеем $((x+2)*3)-(x*3) \equiv \mathrm{const}$.
По предложению~\ref{prozk21} находим, что $x*3 \equiv c\cdot x+d$ для некоторых констант $c$ и $d$.
Если $c\neq 1$, то равенство $c\cdot x+d = x$, т.\,е. $x*3=x*0$, имеет решение,
что противоречит тому, что $*$ является квазигруппой.
Следовательно, $c=1$. Из $0*3 = 3$ мы также находим $d=3$.
Таким образом, $x*3 \equiv x+3$, и из (\ref{eq:23xy}) мы видим, что $f(x,y,2,3)\equiv x+y+2+3$.
Аналогично $f(x,y,z,w)\equiv x+y+z+w$ для всех $z$ и $w$.
\proofend

\section{Сублинейные $n$-квазигруппы}\label{s:Z5Z7}

Будем говорить, что $n$-квазигруппа порядка $k$ {\em сублинейная},
если все ее $2$-ретракты изотопны $Z_k$. Следующие два факта
установлены при помощи компьютерного перебора:

\begin{proposition}
Все сублинейные $4$-квазигруппы порядка $5$ разделимы.
Все сублинейные $3$-квазигруппы порядка $7$ разделимы.
\end{proposition}

Первый факт был проверен непосредственно, используя классификацию
\cite{RU:MK-W:small}, \cite{RU:MK-W:data} $4$"=квазигрупп (латинских $4$-кубов)
порядка $5$. Все сублинейные $3$"=квазигруппы $f$ порядка $7$ были
проверены, используя следующий подход. Сначала мы можем взять, без
потери общности, $f(0,y,z) \equiv y+z$ и $f(x,0,0) \equiv x$. Затем,
имеется $5!=120$ способов выбрать ретракт $f(x,y,0)$, изотопный
$Z_7$, и $120$ способов выбрать $f(x,0,z)$. И наконец, мы выбираем
перестановку $f(1,1,z)$, после чего, по следствию~\ref{cor:z5z5},
сублинейная $3$"=квазигруппа должна восстанавливаться однозначно.

По индукции, с помощью теоремы~\ref{th}, получаем следующее:

\begin{corollary}\label{cor:subli}
Все сублинейные $n$-квазигруппы порядка $5$ разделимы при $n\geq 4$.
Все сублинейные $n$-квазигруппы порядка $7$ разделимы при $n\geq 3$.
\end{corollary}

Имеется единственная, с точностью до изотопии, неразделимая
сублинейная $3$"=квазигруппа порядка $5$. Интересно было бы обобщить
следствие~\ref{cor:subli} на другие порядки. Учитывая
теорему~\ref{th}, в случае простого порядка достаточно показать
разделимость для некоторой фиксированной арности, скажем $3$ или
$4$.

